\documentclass[notitlepage,11pt,a4paper]{scrartcl}
\usepackage{amstext}
\usepackage{amssymb}
\usepackage{amsfonts}
\usepackage{amssymb}
\usepackage{amsmath}
\usepackage{graphicx}
\usepackage[all]{xy}
\usepackage[latin1]{inputenc}
\usepackage{theorem}
\usepackage{geometry}
\usepackage{abstract}
\usepackage{textcomp}
\newtheorem{satz}{Theorem}[section]

\newtheorem{defi}[satz]{Definition}
\newtheorem{bem}[satz]{Remark}
\newtheorem{kor}[satz]{Corollary}
\newtheorem{lem}[satz]{Lemma}

\newtheorem{bei}[satz]{Example}
\newtheorem{pro}[satz]{Proposition}

\newcommand{\qed}{\begin{flushright}$\square$\end{flushright}}
\newcommand{\filt}[6]{
\[
\begin{xy}
\xymatrix@R20pt@C20pt{
&\mathbb{C}^3&\\\langle #1,#2 \rangle\ar[ru]&\langle #3,#4\rangle\ar[u]&\langle #5,#6\rangle\ar[lu]\\\langle #1\rangle\ar[u]&\langle #3\rangle\ar[u]&\langle #5\rangle\ar[u]}
\end{xy}
\]
}

\newcommand\Qn{\mathbb{Q}}
\newcommand\Zn{\mathbb{Z}}
\newcommand\Nn{\mathbb{N}}
\newcommand\Cn{\mathbb{C}}

\newcommand{\Sc}[2]{\langle #1,#2\rangle}
\newcommand{\Hom}{\mathrm{Hom}} 
\newcommand{\Ext}{\mathrm{Ext}}

\newcommand{\PrL}{{}^\prime\!\Lambda}
\newcommand{\ses}[3]{0\rightarrow #1\rightarrow #2\rightarrow#3\rightarrow 0}
\begin{document}
\parindent0pt
\title{Refined GW/Kronecker correspondence}
\author{Markus Reineke and Thorsten Weist\\ Fachbereich C - Mathematik\\ Bergische Universit\"at Wuppertal\\ D - 42097 Wuppertal\\ reineke@math.uni-wuppertal.de, weist@math.uni-wuppertal.de}
\maketitle
\begin{abstract} Gromov-Witten invariants of weighted projective planes and Euler characteristics of moduli spaces of representations of bipartite quivers are related via the tropical vertex, a group of formal automorphisms of a torus. On the Gromov-Witten side, this uses the work of Gross, Pandharipande and Siebert. The quiver moduli side features quiver wall-crossing formulas, functional equations for Euler characteristics, and localization techniques. We derive several explicit formulas for Gromov-Witten invariants.
\end{abstract}
 
\section{Introduction}

It was noted in \cite{GPS}, \cite{RePoisson}, and explained in detail in \cite{GP}, that there is a numerical correspondence (named GW/Kronecker correspondence in \cite{Stoppa}) between two seemingly unrelated geometries :\\[1ex]
The first geometry is the Gromov-Witten theory of weighted projective planes, with Gromov-Witten invariants counting rational maps to open parts of weighted projective planes intersecting the toric divisors in prescribed points with prescribed multiplicities (see \cite[Section 0]{GPS}). The second geometry is the theory of moduli spaces of quiver representations (see for example \cite{ReICRA}), and the relevant numerical information is just the Euler characteristic of these moduli.\\[1ex]
The relation between these geometries is realized in the tropical vertex, a group of formal automorphisms of a torus which first appeared in \cite{KSAffine} and plays a prominent role in describing the wall-crossing behaviour of Donaldson-Thomas invariants in
\cite{KSWall}. Namely, generating series of Gromov-Witten invariants, respectively of Euler characteristics of quiver moduli, appear in a natural factorization of a commutator in the tropical vertex.\\[1ex]
On the Gromov-Witten side, this result of \cite{GPS} requires a passage from the above Gromov-Witten invariants via degeneration formulas and multiple cover calculations in Gromov-Witten theory, and holomorphic and tropical curve counts, to scattering diagrams describing factorizations in the tropical vertex.\\[1ex]
On the quiver moduli side, the wall-crossing formula of \cite{RePoisson} is derived via counting points of quiver moduli over finite fields, realizing the Harder-Narasimhan recursion as an identity in Hall algebras of quivers, and translating this identity to a group of Poisson automorphisms using framed versions of quiver moduli. The tropical vertex appears in the special situation of moduli spaces for generalized Kronecker quivers.\\[2ex]
In this paper, we refine and generalize the GW/Kronecker correspondence, proving that the Euler characteristics of appropriate quiver moduli determine, and are determined by, the above Gromov-Witten invariants. In several cases, this allows us to give explicit formulas for the Gromov-Witten invariants, confirming in particular the conjecture of \cite[Section 1.4]{GP} on Gromov-Witten invariants of the projective plane, and some cases of the integrality conjecture \cite[Conjecture 6.2]{GPS}.\\[1ex]
To achieve this, we first apply the wall-crossing formula of \cite{RePoisson} to complete bipartite quivers; similarly to \cite{GP}, this gives a correspondence between Gromov-Witten invariants and Euler characteristics of framed moduli spaces of representations of these quivers. We then apply the formalism of \cite{ReFE}, relating Euler characteristics of framed and unframed quiver moduli by systems of functional equations for their generating series, resulting in a correspondence between Gromov-Witten invariants and Euler characteristics of unframed moduli spaces of representations of bipartite quivers. Explicit formulas are obtained via localization theory for quiver moduli as developed in \cite{wei}; several of the main results of this paper are generalized to bipartite quivers.\\[2ex]
In sections \ref{tvgw} and \ref{gwi}, we recall the definition of the tropical vertex and the basic factorization problem (\ref{factorization}) which is solved by the two different above geometries, as well as the relevant Gromov-Witten theory. We follow the notation of \cite{GPS} very closely to allow the reader a direct comparison of our results with \cite{GP,GPS}.\\[1ex]
All notions and results from the theory of representations of quivers and their moduli spaces (in particular, the definition and basic geometric properties of the unframed and framed moduli spaces, the wall-crossing formula of \cite{RePoisson} and the functional equations of \cite{ReFE}), which are necessary for the derivation of the refined GW/Kronecker correspondence, are reviewed in Section \ref{qmwcfe}. In sections \ref{cbq} and \ref{rc}, we specialize these methods to complete bipartite quivers, and obtain a first version of the refined GW/Kronecker correspondence in Theorem \ref{gwmk1}. The second version of the correspondence, Theorem \ref{gwmk2} in Section \ref{fe}, is obtained by specializing the functional equations of \cite{ReFE} to bipartite quivers.\\[1ex]
Although this second version shows that Gromov-Witten invariants and Euler characteristics of quiver moduli determine each other, the precise mechanism for this determination is (in general) hidden under an infinite system of coupled functional equations, reminiscent of $Q$-system type equations \cite{KNT} arising from the Bethe ansatz equations of solvable lattice models (although no potential lattice model corresponding to the equations of Theorem \ref{gwmk2} is known). To extract more specific information from Theorem \ref{gwmk2}, special cases where the Euler characteristics of quiver moduli can be computed (usually by localization techniques) are considered in the following sections:\\[1ex]
A trivial first order analysis of the functional equations yields a much simpler correspondence in the so-called coprime case in Section \ref{coprime}: on the Gromov-Witten side, the numbers of intersections with two of the toric divisors of the rational curves to be counted are assumed to be coprime; on the quiver moduli side, coprimality of dimension types of representations results in compact moduli. The correspondence then just states that the two invariants are equal (Corollary \ref{cor1}). For one particular slope (the ratio of the total intersection multiplicities with two toric divisors), localization methods yield an explicit formula for non-trivial counting invariants (Theorem \ref{dm1d}); a conceptual reason for this particular case to be computable (and for being almost the only such example) is still unknown.\\[1ex]
In Section \ref{smalllength}, we consider those (few) cases where the complete bipartite quiver is of (extended) Dynkin type; the known classification of all representations in these cases allows us to give a complete description of the factorization (\ref{factorization}) in the tropical vertex, thereby determining all Gromov-Witten invariants counting curves intersecting the toric divisors in ``few'' points.\\[1ex]
It was already noted in \cite{GP,GPS} via computer experiments that a closed formula, conjectured in \cite[Section 1.4]{GP}, can be expected for the Gromov-Witten invariants of the (unweighted) projective plane. We prove this conjecture in Section \ref{central} using a vanishing result, Corollary \ref{vanish}, for certain Euler characteristics (although the corresponding quiver moduli are highly non-trivial), which simplifies the functional equations to a finite set of algebraic functional equations in Theorem \ref{fecs}. A more detailed closed formula can be expected to follow from an appropriate application of multivariate Lagrange inversion.\\[1ex]
To address integrality properties of the Gromov-Witten invariants and their associated BPS state counts (see \cite[Conjecture 6.2]{GPS}), we consider specialized variables in the tropical vertex, which on the Gromov-Witten side means that only the total numbers of intersections of curves with the toric divisors are recorded. We show in Section \ref{special} that the system of functional equations reduces to a single such equation (Theorem \ref{fes}) of the type which was already studied in \cite{ReFE} in connection with the relative integrality of Donaldson-Thomas invariants. The methods of \cite{ReFE}, together with a subtle divisibility property for the Euler characteristic of quiver moduli (Theorem \ref{divi}) obtained again by a localization argument, yield integrality in the so-called balanced case (see Section \ref{balanced}).\\[1ex]
We briefly consider more general commutator formulas in the tropical vertex \cite[Theorem 5.6]{GPS} in Section \ref{gcf}. This requires consideration of moduli spaces for bipartite quivers with level structure. However, the analogue of the GW/Kronecker correspondence, Theorem \ref{moregeneral}, is weaker in this generality; the reason for this is unclear at the moment.\\[1ex]
Sections \ref{ecms} to \ref{appl} develop the localization techniques which, as indicated above, form the main technical tool for derivation of explicit formulas for Gromov-Witten invariants. After reviewing general concepts in Section \ref{ecms}, the two explicit formulas for Euler characteristics leading to Theorem \ref{dm1d} and Corollary \ref{integ} are developed in Section \ref{Kl1l2}. Several situations (in the generality of levelled bipartite quivers) in which Euler characteristics of quiver moduli vanish are derived in Section \ref{vanisheuler}. Finally, Section \ref{appl} details some further classes of examples.\\[2ex]
In conclusion, the refined GW/Kronecker correspondence, together with localization techniques, allows to obtain several nontrivial formulas for Gromov-Witten invariants. These results hopefully serve as a starting point for investigation of a direct geometric relation between the two geometries in question.\\[3ex]
{\it Acknowledgments:} The authors would like to thank A. King, R. Pandharipande, B. Siebert, Y. Soibelman and J. Stoppa for valuable discussions about the material developed here.

\section{The tropical vertex}\label{tvgw}

We review the definition of the tropical vertex following \cite[Section 0]{GPS}.\\[1ex]
We fix nonnegative integers $l_1,l_2\geq 1$ and define $R$ as the formal power series ring $R=\Qn [[s_1,\ldots,s_{l_1},t_1,\ldots,t_{l_2}]]$, with maximal ideal $\mathfrak{m}$. Let $B$ be the $R$-algebra $$B=\Qn[x^{\pm 1},y^{\pm 1}][[s_1,\ldots,s_{l_1},t_1,\ldots,t_{l_2}]]=\Qn[x^{\pm 1},y^{\pm 1}]\widehat{\otimes}R$$
(a suitable completion of the tensor product). We consider $R$-linear automorphisms of $B$ (more precisely, we only consider automorphisms respecting the symplectic form $\frac{dx}{x}\wedge\frac{dy}{y}$):\\[1ex]
For $(a,b)\in\Zn^2$ and a series $f\in 1+x^ay^b\Qn[x^ay^b]\widehat{\otimes}\mathfrak{m}$, define the automorphism
$$T_{(a,b),f}:\left\{\begin{array}{ccc}x&\mapsto&xf^{-b}\\ y&\mapsto&yf^a.\end{array}\right.$$

\begin{defi} The tropical vertex group $H\subset{\rm Aut}_R(B)$ is defined as the completion with respect to $\mathfrak{m}$ of the subgroup of ${\rm Aut}_R(B)$ generated by all elements $T_{(a,b),f}$ as above.
\end{defi}

By \cite{KSAffine} (see also \cite[Theorem 1.3]{GPS}), there exists a unique factorization in $H$ into an infinite ordered product
\begin{equation}\label{factorization} T_{(1,0),\prod_k(1+s_kx)}T_{(0,1),\prod_l(1+t_ly)}=\prod_{b/a\mbox{ \footnotesize decreasing}}T_{(a,b),f_{(a,b)}},\end{equation}
the product ranging over all coprime pairs $(a,b)\in\Nn^2$. The main problem addressed in \cite{GPS} is to describe the series $f_{(a,b)}$ appearing in this factorization.\\[1ex]
As an example, we list the cases in which the factorization (\ref{factorization}) actually involves only finitely many terms; these are exactly the cases where $l_1l_2\leq 3$. Without loss of generality, we can assume $l_1\leq l_2$; we borrow results from Section \ref{smalllength}:\\[1ex]
If $l_1=1$, $l_2=1$, we have $T_{(1,0),1+s_1x}T_{(0,1),1+t_1y}=$
$$=T_{(0,1),1+t_1y}T_{(1,1),1+s_1t_1xy}T_{(1,0),1+s_1x}.$$
If $l_1=1$, $l_2=2$, we have $T_{(1,0),1+s_1x}T_{(0,1),(1+t_1y)(1+t_2y)}=$
$$=T_{(0,1),(1+t_1y)(1+t_2y)}T_{(1,2),1+s_1t_1t_2xy^2}T_{(1,1),(1+s_1t_1xy)(1+s_1t_2xy)}T_{(1,0),1+s_1x}.$$
If $l_1=1$, $l_2=3$, we have $T_{(1,0),1+s_1x}T_{(0,1),(1+t_1y)(1+t_2y)(1+t_3y)}=$
$$=T_{(0,1),(1+t_1y)(1+t_2y)(1+t_3y)}T_{(1,3),1+s_1t_1t_2t_3xy^3}T_{(1,2),(1+s_1t_1t_2xy^2)(1+s_1t_1t_3xy^2)(1+s_1t_2t_3xy^2)}$$
$$T_{(2,3),1+s_1^2t_1t_2t_3x^2y^3}T_{(1,1),(1+s_1t_1xy)(1+s_1t_2xy)(1+s_1t_3xy)}T_{(1,0),1+s_1x}.$$

There are two more cases where the factorization (\ref{factorization}), although involving infinitely many nontrivial factors, can be described completely; these are the cases $l_1l_2=4$ to be discussed in Section \ref{smalllength}.

\section{Gromov-Witten invariants}\label{gwi}

The main result of \cite{GPS} is the description of the series $f_{(a,b)}$ appearing in the factorization (\ref{factorization}) in terms of Gromov-Witten theory of certain toric surfaces:\\[2ex]
We follow the notation of \cite[Section 0.4]{GPS}. Let $\Sigma\subset\Zn^2$ be the fan with rays generated by $-(1,0),-(0,1),(a,b)$. Let $X_{a,b}$ be the toric surface over $\Cn$ associated to $\Sigma$ (which is isomorphic to the weighted projective plane $(\Cn^3\setminus\{0\})/\Cn^*$ for the action $t(x,y,z)=(t^ax,t^by,tz)$) with corresponding toric divisors $D_1,D_2,D_{\rm out}$. Let $X^o_{a,b}\subset X_{a,b}$ be the open surface obtained by removing the toric fixed points, and let $D^o_1,D^o_2,D^o_{\rm out}$ be the restrictions of the toric divisors to $X^o_{a,b}$.\\[1ex]
We consider a pair $({\bf P}_1,{\bf P}_2)$ of ordered partitions (that is, partitions with nonnegative parts which are allowed to be $0$, and whose order is kept track of), written
$${\bf P}_1=p_{1,1}+\ldots+p_{1,l_1},\;\;\;{\bf P}_2=p_{2,1}+\ldots+p_{2,l_2},$$
such that $|{\bf P}_1|=\sum_{l=1}^{l_1}p_{1,l}=ka$ and $|{\bf P}_2|=\sum_{l=1}^{l_2}p_{2,l}=kb$ for some $k\geq 1$. Let $\nu:X_{a,b}[({\bf P}_1,{\bf P}_2)]\rightarrow X_{a,b}$ be the blow-up of $X_{a,b}$ along $l_1$ (resp.~$l_2$) points of $D_1^o$ (resp. $D_2^o$), and define $X_{a,b}^o[({\bf P}_1,{\bf P}_2)]=\nu^{-1}(X_{a,b}^o)$. Let $\beta_k\in H^2(X_{a,b},\Zn)$ be the unique cohomology class with intersection numbers $$\beta_1\cdot D_1=ka,\; \beta_2\cdot D_2=kb,\; \beta_k\cdot D_{\rm out}=k.$$
Define a cohomology class $\beta_k[({\bf P}_1,{\bf P}_2)]\in H^2(X_{a,b}[({\bf P}_1,{\bf P}_2)],\Zn)$ by
$$\beta_k[({\bf P}_1,{\bf P}_2)]=\nu^*(\beta_k)-\sum_{k=1}^{l_1}p_{1,k}[E_{1,k}]-\sum_{l=1}^{l_2}p_{2,l}[E_{2,l}],$$
where $E_{i,k}$ for $k=1,\ldots,l_i$ denotes the $k$-th exceptional divisor over $D_i^o$ for $i=1,2$. The moduli space $\overline{\mathfrak{M}}(X_{a,b}^o[({\bf P}_1,{\bf P}_2)]/D^o_{\rm out})$ of genus $0$ maps to $X_{a,b}^o[({\bf P}_1,{\bf P}_2)]$ in class $\beta_k[({\bf P}_1,{\bf P}_2)]$ with full contact order $k$ at an unspecified point of $D_{\rm out}^o$ is proper and of virtual dimension $0$, thus a corresponding Gromov-Witten invariant $N_{a,b}[({\bf P}_1,{\bf P}_2)]\in\Qn$ is well-defined (see \cite[Section 5.2]{GPS}).\\[1ex]
Heuristically, $N_{a,b}[({\bf P}_1,{\bf P}_2)]$ may be viewed as the ``number''  of rational curves in $X_{a,b}$ intersecting the distinct fixed $l_i$ points of $D_i^o$ with multiplicities given by the $p_{i,l}$ for $i=1,2$, and being tangent to $D_{\rm out}^o$ of order $k$. But note that this counting is only straightforward in very particular cases; in general, complicated contributions from degenerations and multiple covers have to be accounted for.

\begin{bei}
\end{bei}
We have $N_{(1,3)}[(1,1+1+1)]=1$, since there is (up to reparametrization) a unique curve
$$(u:v)\mapsto(u:-\frac{y_1}{x_1x_2x_3}(u-x_1v)(u-x_2v)(u-x_3v):v)$$
from ${\bf P}^1$ to $X_{1,3}$ intersecting $D^o_1$ transversally in three points $x_1,x_2,x_3$, intersecting $D^o_2$ transversally in one point $y_1$, and intersecting $D^o_{\rm out}$ transversally in one point.\\[1ex]
We also have $N_{(2,3)}[2,1+1+1]=1$, but the counting procedure is more indirect: choosing signs $\varepsilon_1,\varepsilon_2,\varepsilon_3$ such that $\varepsilon_1\varepsilon_2\varepsilon_3=1$, there are four curves
$$(u:v)\mapsto(u^2:-\frac{y_1}{x_1x_2x_3}(u-\varepsilon_1x_1v)(u-\varepsilon_2x_2v)(u-\varepsilon_3x_3v):v)$$
from ${\bf P}^1$ to $X_{2,3}$ intersecting $D^o_1$ transversally in three points $x_1^2,x_2^2,x_3^2$, intersecting $D^o_2$ with multiplicity two in one point $y_1$, and intersecting $D^o_{\rm out}$ transversally in one point. By the degeneration formula \cite[Proposition 5.3]{GPS}, these have to be weighted by the factor $-\frac{1}{2}$ in the computation of $N_{(2,3)}[2,1+1+1]$, and six curves intersecting $D^o_2$ transversally in two distinct points have to be added with a weight factor of $\frac{1}{2}$.\\[1ex]
We have $N_{(1,1)}[(1+1,1+1)]=2$, since there are two curves
$$(u:v)\mapsto(u(u-v):(u-2v)(u-4v):v^2),$$
$$(u:v)\mapsto(u(u-\frac{5}{\sqrt{3}}v):-(u-2\sqrt{3}v)(u+\frac{4}{\sqrt{3}}v):v^2)$$
intersecting $D^o_1$ transversally in the two points $2$, $12$, intersecting $D^o_2$ transversally in the two points $3$, $8$, and being tangent of order $2$ to $D^o_{\rm out}$.\\[1ex]
See also \cite[Section 6.4]{GPS} for some more examples.\\[3ex]
The main result of \cite{GPS} is

\begin{satz}\cite[Theorem 5.4]{GPS}\label{t54} For all coprime $(a,b)$, we have
$$\log f_{(a,b)}=\sum_{k\geq 1}\sum_{\substack{{|{\bf P}_1|=ka,}\\{|{\bf P}_2|=kb}}}kN_{(a,b)}[({\bf P}_1,{\bf P}_2)]s^{{\bf P}_1}t^{{\bf P}_2}(x^ay^b)^k.$$
\end{satz}

\section{Recollections on quiver moduli, wall-crossing and functional equations}\label{qmwcfe}

Let $Q$ be a finite quiver, given by a finite set $Q_0$ of vertices and finitely many arrows $\alpha:i\rightarrow j$. We denote by $\Lambda=\Zn Q_0$ the free abelian group over $Q_0$, by $\Lambda^+=\Nn Q_0\subset\Lambda$ the set of dimension vectors, which will be written as $d=\sum_{i\in Q_0}d_ii\in\Lambda^+$, and by $\PrL^+=\Lambda^+\setminus\{0\}$ the set of nonzero dimension vectors. The Euler form on $\Lambda$ is given by $$\langle d,e\rangle=\sum_{i\in Q_0}d_ie_i-\sum_{\alpha:i\rightarrow j}d_ie_j.$$
A representation $V$ of $Q$ of dimension vector $d\in\Lambda$ (over the field of complex numbers) consists of complex vector spaces $V_i$ for $i\in Q_0$ of dimension $d_i$, and of linear maps $V_\alpha:V_i\rightarrow V_j$ for every arrow $\alpha:i\rightarrow j$ in $Q$.\\[1ex]
We discuss the roots of a quiver; see \cite{kac} for more details. A dimension vector is called a root if there exists at least one indecomposable representation of this dimension and it is called Schur root if, in addition, the endomorphism ring of at least one representation is trivial. The last condition already implies that an open subset of representations has trivial endomorphism ring, see for instance \cite{sch}. A root $d\in\Nn Q_0$ is called real if we have $d\in W(Q)Q_0$, i.e. $d$ is a reflection of a simple root, where $W(Q)$ denotes the Weyl group of the quiver. All the other roots are called imaginary. It is well-known that a root is real if and only if $\Sc{d}{d}=1$ and imaginary if and only if $\Sc{d}{d}\leq 0$. If $\Sc{d}{d}=0$, we also call $d$ isotropic.\\[1ex]
Let $(d,e):=\Sc{d}{e}+\Sc{e}{d}$ be the symmetrized Euler form. The fundamental domain $F(Q)$ of $\mathbb{N}Q_0$ is given by the dimension vectors $d$ such that $(d,q)\leq 0$ for all $q\in Q_0$. Moreover, we have $d\in W(Q)F(Q)$ for all imaginary roots $d$.\\[1ex]
Let $\Theta\in\Lambda^*$ be a functional on $\Lambda$, viewed as a stability. It induces a slope function $\mu:\PrL^+\rightarrow\Qn$ by $\mu(d)=\Theta(d)/\dim d$, where $\dim\in\Lambda^*$ is the functional $\dim d=\sum_{i\in Q_0}d_i$ (in Section \ref{gcf}, we will consider more general stabilities: additionally to the functional $\Theta$, we choose a functional $\kappa$ assuming positive values on $\PrL^+$ and define the slope function by $\mu(d)=\Theta(d)/\kappa(d)$; all results reviewed in the following are easily seen to generalize to this context with obvious modifications). For $\mu\in\Qn$, let $\PrL_\mu^+$ be the set of dimension vectors $d\in\PrL^+$ of slope $\mu$, and let $\Lambda^+_\mu$ be $\PrL_\mu^+\cup\{0\}$, a subsemigroup of $\Lambda^+$.\\[1ex]
The slope of a representation $V$ of $Q$ is defined as the slope of its dimension vector. The representation $V$ is called $\Theta$-(semi-)stable if the slope (weakly) decreases on proper non-zero subrepresentations.\\[1ex]
There exists a moduli space $M_d^{\Theta-{\rm st}}(Q)$ for isomorphism classes of $\Theta$-stable representations of $Q$ of dimension vector $d$. If non-empty, it is a smooth irreducible variety of dimension $1-\langle d,d\rangle$. It is projective if $Q$ has no oriented cycles and $d$ is $\Theta$-coprime, that is, if $\mu(e)=\mu(d)$ for $0\not=e\leq d$ implies $e=d$. Call a dimension vector $d\in\Lambda^+$ indivisible if ${\rm gcd}(d_i\, :\, i\in Q_0)=1$. For generic $\Theta$, a dimension vector $d$ is $\Theta$-coprime if and only if it is indivisible.\\[1ex]
For a $\Theta$-coprime dimension vector $d$, there is an explicit formula for the Poincar\'e polynomial of $M_d^{\Theta-{\rm st}}(Q)$ arizing from a resolution of a Harder-Narasimhan type recursion:

\begin{satz}\label{rhnr}\cite[Corollary 6.8]{ReHNS} For $\Theta$-coprime $d$, we have
$$\sum_i\dim H^i(M_d^{\Theta-{\rm st}}(Q),\Qn)q^{i/2}=(q-1)\sum_{d^*}(-1)^{s-1}q^{-\sum_{k\leq l}\langle d^l,d^k\rangle}\prod_{k=1}^s\prod_{i\in Q_0}\prod_{j=1}^{d^k_i}(1-q^{-j})^{-1},$$
where the sum ranges over all decompositions $d=d^1+\ldots+d^s$ of $d$ such that all $d^k$ are non-zero, and $\mu(d^1+\ldots+d^k)>\mu(d)$ for all $k<s$.
\end{satz}

In sections \ref{ecms}, \ref{Kl1l2}, \ref{vanisheuler}, we will make key use of localization theory for quiver moduli as developed in \cite{wei}. The starting point for this method is the fact that the Euler characteristic $\chi(M_d^{\Theta-st}(Q))$ for a quiver containing (unoriented) cycles can be computed in terms of the Euler characteristics of various moduli spaces of representations of the universal covering $\tilde{Q}$ of $Q$ (for precise definitions and the proof, see Section \ref{ecms}):

\begin{satz}\label{euler} Let $Q$ be a quiver with dimension vector $d$. Then for the Euler characteristic of the moduli space $M^{\Theta-\rm{st}}_d(Q)$ we have
\[\chi(M^{\Theta-\rm{st}}_d(Q))=\sum_{\tilde{d}}\chi(M^{\tilde{\Theta}-\rm{st}}_{\tilde{d}}(\tilde{Q})),\]
where $\tilde{d}$ ranges over all equivalence classes being compatible with $d$.\end{satz}
Let $n\in\Lambda^+$ be another dimension vector, and choose a complex vector space $W_i$ of dimension $n_i$ for each vertex $i\in Q_0$. There exists a moduli space $M_{d,n}^\Theta(Q)$ parametrizing equivalence classes of pairs $(V,f)$, where $V$ is a $\Theta$-semistable representation of $Q$ of dimension vector $d$ and $f=(f_i:W_i\rightarrow V_i)_{i\in Q_0}$ is an $Q_0$-graded linear map such that the following holds: if $U\subset V$ is a subrepresentation containing the image of $f$, that is, $f_i(W_i)\subset U_i$ for all $i\in Q_0$, then $\mu(U)<\mu(V)$. Such objects are parametrized up to isomorphisms of semistable representations intertwining the $Q_0$-graded linear maps. This moduli space is called a smooth model in \cite{ERHilb2}. If non-empty, $M_{d,n}^\Theta(Q)$ is a smooth irreducible variety of dimension $n\cdot d-\langle d,d\rangle$, where $n\cdot d=\sum_{i\in Q_0}n_id_i$. It contains an open subset which is a ${\bf P}^{n\cdot d-1}$-fibration over $M_d^{\Theta-{\rm st}}(Q)$; thus it can be viewed as a (partial) compactification of a projective space fibration over $M_d^{\Theta-{\rm st}}(Q)$ -- note that no natural smooth compactifications of $M_d^{\Theta-{\rm st}}(Q)$ are known in general.\\[1ex]
We define the generating series of Euler characteristics
$$Q_\mu^{(n)}(x)=\sum_{d\in\Lambda^+_\mu}\chi(M_{d,n}^\Theta(Q))x^d\in\Zn[[\Lambda_\mu^+]].$$
Let $\{d,e\}=\langle d,e\rangle-\langle e,d\rangle$ be the antisymmetrization of the Euler form. We assume that $Q$ has no oriented cycles, thus we can order the vertices as $Q_0=\{i_1,\ldots,i_r\}$ in such a way that $k > l$ provided there exists an arrow $i_k\rightarrow i_l$.\\[1ex]
Define a Poisson algebra $B(Q)=\Qn[[x_i\, :\, i\in Q_0]]$ with Poisson bracket $\{x^d,x^e\}=\{d,e\}x^{d+e}$, where $x^d=\prod_{i\in Q_0}x_i^{d_i}$ denotes the natural topological basis of $B(Q)$. We consider Poisson automorphisms of $B(Q)$. For a vertex $i\in Q_0$, define $T_{i}\in{\rm Aut}(B(Q))$ by $T_i(x^d)=x^d(1+x_i)^{\{i,d\}}$. 

\begin{satz}\label{wcqm}\cite[Theorem 2.1]{RePoisson} We have the following factorization in ${\rm Aut}(B(Q))$:
$$T_{i_1}\circ\ldots\circ T_{i_r}=\prod_{\mu\in\Qn}^\leftarrow T_\mu,$$
where
$$T_\mu(x^d)=x^d\prod_{i\in Q_0}Q_\mu^i(x)^{\{i,d\}}.$$
\end{satz}

For an arbitrary functional $\eta\in(\Qn Q_0)^*$, we define $Q_\mu^{\eta}(x)=\prod_{i\in Q_0}Q_\mu^i(x)^{\eta(i)}$, thus in particular $Q_\mu^{n\cdot}(x)=Q_\mu^{(n)}(x)$. By \cite[Lemma 3.6]{RePoisson}, we have:

\begin{lem}\label{identity} The series $Q_\mu^{\Theta-\mu\dim}$ equals $1$.
\end{lem}

The following is the main result of \cite{ReFE}:

\begin{satz}\label{fneq} The series $Q_\mu^{(n)}(x)$ is given by
$$Q_\mu^{(n)}(x)=\prod_{d\in'\Lambda^+_\mu}R^d(x)^{\chi(d)(n\cdot d)},$$
where the series $R^d(x)\in\Zn[[\Lambda_\mu^+]]$ for $d\in\PrL_\mu^+$ are uniquely determined by the following system of functional equations:\\[1ex]
For all $d\in\PrL_\mu^+$, we have
$$R^d(x)=(1-x^d\prod_{e\in'\Lambda^+_\mu}R^e(x)^{-\chi(e)\langle d,e\rangle})^{-1},$$
where $\chi(d)=\chi(M_d^{\Theta-\rm st}(Q))$.
\end{satz}

\section{Complete bipartite quivers}\label{cbq}

In this section, we specialize the results of the previous section to a class of complete bipartite quivers.\\[1ex]
Let $K=K(l_1,l_2)$ be the quiver with set of vertices $Q_0=\{ i_1,\ldots,i_{l_1},j_1,\ldots,j_{l_2}\}$, and one arrow from each vertex $j_l$ to each vertex $i_k$; thus $K$ is a complete bipartite quiver. We thus have a natural $S_{l_1}\times S_{l_2}$-symmetry of $K(l_1,l_2)$. Dimension vectors for $K$ can be written as $d({\bf P}_1,{\bf P}_2)$ for a pair $({\bf P}_1,{\bf P}_2)$ of ordered partitions ${\bf P}_i=p_{i,1}+\ldots+p_{i,l_i}$ of length $l_i$ (for $i=1,2$) as in Section \ref{gwi}, via $d({\bf P}_1,{\bf P}_2)_{i_k}=p_{1,k}$, $d({\bf P}_1,{\bf P}_2)_{j_l}=p_{2,l}$. We will always assume that $|{\bf P}_1|=ka$, $|{\bf P}_2|=kb$ for $a$ and $b$ coprime and $k\geq 1$. The Euler form on $K(l_1,l_2)$ is given by $$\langle d({\bf P}_1,{\bf P}_2),d({\bf P}_1',{\bf P}_2')\rangle=\sum_k p_{1,k}p_{1,k}'+\sum_l p_{2,l}p_{2,l}'-\sum_{k,l}p_{1,k}'p_{2,l},$$
and thus its antisymmetrization is given by
$$\{d({\bf P}_1,{\bf P}_2),d({\bf P}_1',{\bf P}_2')\}=\sum_{k,l}(p_{1,k}p_{2,l}'-p_{1,k}'p_{2,l}).$$
We choose the stability $\Theta$ on $K(l_1,l_2)$ given by $\Theta(j_l)=1$, $\Theta(i_k)=0$ (in fact, among the stabilities respecting the symmetry of the quiver, this is the only nontrivial one, see \cite[Section 5.1]{ReICRA} for a discussion in the case of the $m$-Kronecker quiver $K(m)$).\\[1ex]
We describe the resulting moduli spaces:\\[1ex]
Let $V_k$ be $\Cn$-vector spaces of dimension $p_{1,k}$ for $k=1,\ldots,l_1$, and let $W_l$ be vector spaces of dimension $p_{2,l}$ for $l=1,\ldots,l_2$, respectively. Then we consider the action of the group $G=\prod_{k=1}^{l_1}{\rm GL}(V_k)\times\prod_{l=1}^{l_2}{\rm GL}(W_l)$ on the space $\bigoplus_{k,l}{\rm Hom}(W_l,V_k)$ of tuples $(f_{k,l}:W_l\rightarrow V_k)_{k,l}$ of linear maps by base change. Such a tuple is called stable if for all proper non-zero tuples $(U_l\subset W_l)_l$ of subspaces, we have
$$\sum_k\dim\sum_l f_{kl}(U_l)>\frac{a}{b}\sum_l\dim U_l.$$

\begin{satz} The moduli space $$M^{\rm st}({\bf P}_1,{\bf P}_2)=M^{\Theta-{\rm st}}_{d({\bf P}_1,{\bf P}_2)}(K(l_1,l_2))$$ parametrizes stable tuples $(f_{k,l})_{k,l}$ up to the action of $G$. If non-empty, it is a smooth and irreducible variety of dimension $$1-\sum_i p_{1,i}^2-\sum_j p_{2,j}^2+\sum_{k,l}p_{1,k}p_{2,l}.$$
It is projective if $k=1$, that is, if $|{\bf P}_1|$ and $|{\bf P}_2|$ are coprime.
\end{satz}

We discuss the root system of the bipartite quiver $K(l_1,l_2)$. It is easy to check that the inequalities defining the fundamental domain are given by $2p_{1,k}\leq\sum_{l}p_{2,l}$ for all $k$ and $2p_{2,l}\leq \sum_{k}p_{1,k}$ for all $l$.\\[1ex]
Since, in particular, we have $d({\bf P}_1,{\bf P}_2)\notin F(K(l_1,l_2))$ for some real root $d({\bf P}_1,{\bf P}_2)$, we obtain that at least one inequality $2p_{1,k}>\sum_lp_{2,l}$ or $2p_{2,l}>\sum_kp_{1,k}$ holds. Moreover, it is well known that for every real root there exists a (up to isomorphism) unique indecomposable representation.\\[1ex]
But even in the case $l_1=1$, i.e. the $l_2$-subspace quiver, it is not easy to decide if a root is a Schur root. One possibility is to consider the canonical decomposition of the root, see \cite{sch} for more details. If it consists of the dimension vector itself, the dimension vector is a Schur root. Nevertheless, it is possible to recursively construct plenty of real Schur roots of $K(l_1,l_2)$. Therefore, given a real Schur root $d({\bf P}_1,{\bf P}_2)$ of $K(l_1,l_2)$ we consider the quiver $K(l_1,l_2+1)$. Then we have the following lemma:
\begin{lem}
The dimension vector $d({\bf P}_1,\hat{\bf P}_2)$ with $\hat{\bf P}_2=p_{2,1}+\ldots+p_{2,l_2}+(\sum_kp_{1,k})$ is a real Schur root for $K(l_1,l_2+1)$.
\end{lem}
{\it Proof.} It is easy to check that $\Sc{d({\bf P}_1,\hat{\bf P}_2)}{d({\bf P}_1,\hat{\bf P}_2)}=1$. Moreover, the canonical decomposition of $d({\bf P}_1,\hat{\bf P}_2)$ is trivial because the one of $d({\bf P}_1,{\bf P}_2)$ is trivial.
\qed
Note that the setup in the proof is equivalent to the one of the generalized Kronecker quiver $K(\sum_kp_{1,k})$ with dimension vector $(1,\sum_kp_{1,k})$ which also is a real Schur root. Moreover, the preceding statement can be easily generalized to the case of general bipartite quivers.\\[2ex]
For very small cases of $l_1$ and $l_2$, we have a complete description of all moduli spaces. More precisely, for $l_1l_2\leq 4$, the quiver $K(l_1,l_2)$ is of (extended) Dynkin type, and a classification of all isomorphism classes of representations is known:\\[1ex]
If $l_1l_2\leq 3$, the group $G$ has only finitely many orbits even in the space of all tuples of linear maps (the quiver $K(l_1,l_2)$ is then of Dynkin type $A_2$ or $A_3$ or $D_4$). If $l_1l_2=4$, the quiver $K(l_1,l_2)$ is an extended Dynkin quiver of type $\widetilde{A}_3$ or $\widetilde{D}_4$, and again, a classifications of all orbits is known. We use the following lemma which reduces the determination of all moduli spaces to the classification of indecomposable representations with trivial endomorphism ring:

\begin{lem} We have $M^{\rm st}({\bf P}_1,{\bf P}_2)\not=\emptyset$ if and only if there exists a representation of $K(l_1,l_2)$ of dimension vector $d({\bf P}_1,{\bf P}_2)$ with trivial endomorphism ring.
\end{lem}
{\it Proof.} By \cite[Theorem 6.1]{sch}, there exists a representation with trivial endomorphism ring of dimension vector $d=d({\bf P}_1,{\bf P}_2)$ if and only if there exists a stable representation with respect to the stability $\Theta_d=\{d,\-\}$. But the conditions $\mu(e)<\mu(d)$, defined with respect to $\Theta_d$ and $\Theta$, respectively, are equivalent.
\qed

We can now apply the known representation theory of (extended) Dynkin quivers, and in particular the classification of representations with trivial endomorphism ring, to obtain the following description of pairs of ordered partitions (up to reordering) with nonempty moduli spaces:\\[1ex]
$(l_1=1,l_2=1):$
$$ (1,0),\; (1,1),\; (0,1).$$
$(l_1=1,l_2=2):$
$$ (1,0+0),\; (1,1+0),\; (1,1+1),\; (0,1+0).$$
$(l_1=1,l_2=3):$
$$ (1,0+0+0),\; (1,1+0+0),\; (2,1+1+1),\; (1,1+1+0),\; (1,1+1+1),\; (0,1+0+0).$$
$(l_1=1,l_2=4):$
$$ (2i,i+i+i+(i+1)),\; (2i+1,(i+1)+(i+1)+(i+1)+(i+1)),\; (2i+1,i+(i+1)+(i+1)+(i+1)),$$
$$(1,1+1+0+0),\; (2,1+1+1+1),$$
$$(2i+2,i+(i+1)+(i+1)+(i+1)),\; (2i+1,i+i+i+i),\; (2i+1,i+i+i+(i+1)).$$
$(l_1=2,l_2=2):$
$$ (i+(i+1),i+i),\; (i+i,i+(i+1)),$$
$$(0+1,0+1),\; (1+1,1+1),$$
$$((i+1)+(i+1),i+(i+1)),\; (i+(i+1),(i+1)+(i+1)), \mbox{ for }i\geq 0.$$

In almost all of the above cases, the moduli space $M^{\rm st}({\bf P}_1,{\bf P}_2)$ reduces to a single point, except the following:
$$M^{\rm st}(2,1+1+1+1)\simeq{\bf P}^1\setminus\{0,1,\infty\},$$
$$M^{\rm st}(1+1,1+1)\simeq{\bf P}^1\setminus\{0,1\}.$$

We now consider a particular case of the smooth models of the previous section.\\[1ex]
Call a tuple $((f_{k,l}),w_l)$ consisting of a system of linear maps $(f_{k,l}:W_l\rightarrow V_k)_{k,l}$ as above and a system of vectors $(w_l\in W_l)_l$ stable if for all tuples $U_l\subset W_l$ of subspaces, we have $\sum_k\dim\sum_l f_{kl}(U_l)\geq\frac{a}{b}\sum_l\dim U_l$, with strict inequality if $w_l\in U_l$ for all $l=1,\ldots,l_2$. Such tuples are considered up to simultaneous base change in all $V_k$ and all $W_l$, i.e. tuples $(f_{k,l},w_l)$ and $(f_{k,l}',w_l')$ are equivalent if there exist automorphisms $g_k\in{\rm GL}(V_k)$ and $h_l\in{\rm GL}(W_l)$ such that $f_{k,l}'=g_kf_{k,l}h_l^{-1}$ and $w_l'=h_lw_l$ for all $k,l$.\\[1ex]
We define a dimension vector $n^{\rm b}$ for $K(l_1,l_2)$ by $n^{\rm b}_{i_k}=0$ for all $k=1,\ldots,l_1$ and $n^{\rm b}_{j_l}=1$ for all $l=1,\ldots,l_2$.

\begin{satz} The moduli space
$${M}^{\rm b}({\bf P}_1,{\bf P}_2)=M_{d({\bf P}_1,{\bf P}_2),n^{\rm b}}^\Theta(K(l_1,l_2))$$
parametrizes equivalence classes of stable tuples $((f_{k,l}),w_k)$ as above.
\end{satz}

We also have a dual version: we define $n^{\rm f}$ by $n^{\rm f}_{i_k}=1$ for all $k=1,\ldots,l_1$ and $n^{\rm f}_{j_l}=0$ for all $l=1,\ldots,l_2$.
The corresponding moduli space $${M}^{\rm f}({\bf P}_1,{\bf P}_2)=M_{d({\bf P}_1,{\bf P}_2),n^{\rm f}}^\Theta(K(l_1,l_2))$$ parametrizes stable tuples $((f_{k,l}),v_k)$ consisting of a system of linear maps $(f_{k,l}:W_l\rightarrow V_k)_{k,l}$ as above and a system of vectors $(v_k\in V_k)_k$, where stable means that for all tuples $U_l\subset W_l$ of subspaces, we have $\sum_k\dim\sum_l f_{kl}(U_l)\geq\frac{a}{b}\sum_l\dim U_l$, with strict inequality if $v_k\in \sum_{l}f_{kl}(U_l)$ for all $k=1,\ldots,l_1$.\\[1ex]
We define the generating series of Euler characteristics of these moduli spaces by a slight variant of the series $Q_\mu^{(n)}(x)$ of the previous section, namely
$$F^{\rm b}_{(a,b)}=\sum_{k\geq 0}\sum_{\substack{{({\bf P}_1,{\bf P}_2):}\\ {|{\bf P}_1|=ka,|{\bf P}_2|=kb}}}\chi({M}^{\rm b}({\bf P}_1,{\bf P}_2))s^{{\bf P}_1}t^{{\bf P}_2}x^{ka}y^{kb}\in B$$
and
$$F^{\rm f}_{(a,b)}=\sum_{k\geq 0}\sum_{\substack{{({\bf P}_1,{\bf P}_2):}\\ {|{\bf P}_1|=ka,|{\bf P}_2|=kb}}}\chi({M}^{\rm f}({\bf P}_1,{\bf P}_2))s^{{\bf P}_1}t^{{\bf P}_2}x^{ka}y^{kb}\in B$$

\section{Refined GW/Kronecker correspondence}\label{rc}

The main result of this section is the following description of the series $f_{(a,b)}$ appearing in the factorization (\ref{factorization}) in terms of Euler characteristics of the smooth models of complete bipartite quivers of the previous section:

\begin{satz}\label{gwmk1} For all coprime $(a,b)$, we have
$$f_{(a,b)}=(F^{\rm b}_{(a,b)})^{1/b}=(F^{\rm f}_{(a,b)})^{1/a}.$$
\end{satz}
{\it Proof.}
We consider the Poisson algebra $\overline{B}=B(K(l_1,l_2))=\Qn[[x_1,\ldots,x_{l_1},y_1,\ldots,y_{l_2}]]$. It embeds into the algebra $B=\Qn [x^{\pm 1},y^{\pm 1}][[s_1,\ldots,s_{l_1},t_1,\ldots,t_{l_2}]]$ of Section \ref{tvgw} by substituting $x_k=s_kx$ and $y_l=t_ly$.\\[1ex]
The Poisson automorphism $T_{i_k}$ of $\overline{B}$ is then given by 
$$T_{i_k}(x_l)=x_l,\;\;\; T_{i_k}(y_l)=y_l(1+x_k),$$
and similarly $T_{j_k}$ is given by
$$T_{j_k}(x_l)=x_l(1+y_k)^{-1},\;\;\; T_{j_k}(y_l)=y_l.$$
Under the above substitution, this induces automorphisms of $B$ given by
$$T_{i_k}(x)=x,\;\;\; T_{i_k}(y)=y(1+s_kx),$$
$$T_{j_k}(x)=x(1+t_ky)^{-1},\;\;\; T_{j_k}(y)=y.$$
It follows that $T_i=\prod_{k=1}^{l_1}T_{i_k}$ and $T_j=\prod_{l=1}^{l_2}T_{j_l}$ are precisely the automorphisms $T_{(1,0),\prod_{k=1}^{l_1}(1+s_kx)}$ and $T_{(0,1),\prod_{l=1}^{l_2}(1+t_ly)}$, respectively, of Section \ref{tvgw}.\\[1ex]
Now Theorem \ref{wcqm} gives a factorization
$$T_iT_j=\prod_{\mu\in\Qn }^{\leftarrow}T_\mu,$$
where $T_\mu$ is given as follows in terms of the generating functions $Q_\mu^{i_k}(x,y), Q_\mu^{j_l}(x,y)$:
$$T_\mu(x_{i_k})=x_{i_k}\underbrace{\prod_{k'}Q_\mu^{i_{k'}}(x,y)^{\{i_{k'},i_k\}}}_{=1}\prod_lQ_\mu^{j_l}(x,y)^{\{j_l,i_k\}},$$
$$T_\mu(x_{j_l})=x_{j_l}\prod_kQ_\mu^{i_k}(x,y)^{\{i_k,j_l\}}\underbrace{\prod_{l'}Q_\mu^{j_{l'}}(x,y)^{\{j_{l'},j_l\}}}_{=1},$$
the simplifications resulting from the calculation of $\{\_,\_\}$ above.
It follows that $T_\mu$ induces the following automorphism of $B$:
$$T_\mu(x)=x\prod_lQ_\mu^{j_l}(x,y),\;\;\; T_\mu(y)=y\prod_kQ_\mu^{i_k}(x,y).$$
We write $\mu=b/(a+b)$ for $a,b\in\Nn$ coprime. In the series $Q_\mu^{i_k}(x,y)$, $Q_\mu^{j_l}(x,y)$, only monomials $x^{{\bf P}_1}y^{{\bf P}_2}$ with $|{\bf P}_1|=ka$, $|{\bf P}_2|=kb$ for some $k\geq 0$ appear (or, under the above substitution, only monomials $s^{{\bf P}_1}t^{{\bf P}_2}x^{ka}y^{kb}$ appear). Using Lemma \ref{identity}, we have $$(\prod_lQ_\mu^{j_l}(x,y))^{a}(\prod_kQ_\mu^{i_k}(x,y))^{-b}=1.$$
Since $a,b$ are coprime, we can choose $c,d\in\Zn$ such that $ac+bd=1$. We can then define $$G_\mu(x,y)=(\prod_kQ_\mu^{i_k}(x,y))^{c}(\prod_lQ_\mu^{j_l}(x,y))^{d}.$$
A short calculation shows that
$$G_\mu(x,y)^a=\prod_kQ_\mu^{i_k}(x,y),\;\;\; G_\mu(x,y)^b=\prod_lQ_\mu^{j_l}(x,y).$$
Thus, we have written the automorphism $T_\mu$ in the form
$$T_\mu(x)=xG_\mu(x,y)^{-b},\;\;\; T_\mu(y)=yG_\mu(x,y)^{a},$$
i.e. $T_\mu=T_{(a,b),G_\mu(x,y)}$ is an element of the tropical vertex group $H$, and $G_\mu(x,y)$ equals the series $f_{(a,b)}$ by uniqueness of the factorization (\ref{factorization}).\\[1ex]
We have $$G_\mu(x,y)^b=Q_\mu^{(n_0)}(x,y)=F^{\rm b}_{(a,b)}$$
(and similarly for $F^{\rm f}_{(a,b)}$), proving the theorem.
\qed

Comparison of this theorem with Theorem \ref{t54} yields the first instance of the refined GW/Kronecker correspondence:

\begin{kor} For all coprime $(a,b)$, we have

$$\exp\left(\sum_{k=1}^\infty\sum_{{|{\bf P}_1|=ka,|{\bf P}_2|=kb}}kN_{(a,b)}[({\bf P}_1,{\bf P}_2)]s^{{\bf P}_1}t^{{\bf P}_2}x^{ka}y^{kb}\right)=$$
$$=\left(\sum_{k=0}^\infty\sum_{{|{\bf P}_1|=ka,|{\bf P}_2|=kb}}\chi({M}^{\rm b}({\bf P}_1,{\bf P}_2))s^{{\bf P}_1}t^{{\bf P}_2}x^{ka}y^{kb}\right)^{1/b}=$$
$$=\left(\sum_{k=0}^\infty\sum_{{|{\bf P}_1|=ka,|{\bf P}_2|=kb}}\chi({M}^{\rm f}({\bf P}_1,{\bf P}_2))s^{{\bf P}_1}t^{{\bf P}_2}x^{ka}y^{kb}\right)^{1/a}.$$
\end{kor}

Comparing coefficients, we see that all Gromov-Witten invariants $N_{(a,b)}[({\bf P}_1,{\bf P}_2)]$ are determined by all Euler characteristics $\chi(M^{\rm b}({\bf P}_1,{\bf P}_2))$ (or $\chi(M^{\rm f}({\bf P}_1,{\bf P}_2))$) and vice versa.

\section{Derivation of functional equations}\label{fe}

In this section, we apply Theorem \ref{fneq} to derive functional equations determining the series $f_{(a,b)}$.\\[1ex]
For a pair of ordered partitions $({\bf P}_1,{\bf P}_2)$, we abbreviate by $\chi({\bf P}_1,{\bf P}_2)$ the Euler characteristic $\chi(M^{\rm st}({\bf P}_1,{\bf P}_2))$.

\begin{satz}\label{gwmk2} For coprime $(a,b)$, the series $f_{(a,b)}$ is given by
$$f_{(a,b)}=\prod_{k\geq 1}\prod_{\substack{{|{\bf P}_1|=ka}\\{|{\bf P}_2|=kb}}}(R^{{\bf P}_1,{\bf P}_2})^{k\chi({{\bf P}_1,{\bf P}_2})},$$
where the series $R^{{\bf P}_1,{\bf P}_2}\in B$ are determined by the following system of functional equations:\\
For all pairs of ordered partitions $({\bf P}_1,{\bf P}_2)$ as above, 
$$R^{{\bf P}_1,{\bf P}_2}=(1-s^{{\bf P}_1}t^{{\bf P}_2}(x^ay^b)^k\prod_{k'\geq 1}\prod_{\substack{{|{\bf P}_1'|=k'a}\\{|{\bf P}_2'|=k'b}}}(R^{{\bf P}_1',{\bf P}_2'})^{-\langle d({\bf P}_1,{\bf P}_2),d({\bf P}_1',{\bf P}_2')\rangle \chi({\bf P}_1',{\bf P}_2')})^{-1}.$$
\end{satz}
{\it Proof.} In the notation of the proof of Theorem \ref{gwmk1}, we have
$$f_{(a,b)}=G_\mu(x,y)=\prod_kQ_\mu^{i_k}(x,y))^{c}\prod_lQ_\mu^{j_l}(x,y))^{d}=Q_\mu^{n(c,d)\cdot}(x,y),$$
where $n(c,d)\cdot$ is the functional given by $n(c,d)\cdot i_k=c$ for $k=1,\ldots,l_1$ and $n(c,d)\cdot j_l=d$ for $l=1,\ldots,l_2$.
We then have
$$n(c,d)\cdot d({\bf P}_1,{\bf P}_2)=c|{\bf P}_1|+d|{\bf P}_2|=kac+kbd=k.$$
Now the statement of the theorem is just an adaption of Theorem \ref{fneq} to the present notation.
\qed

This second instance of the refined GW/Kronecker correspondence shows (by comparing coefficients $(x^ay^b)^k$) that all Gromov-Witten invariants $N_{(a,b)}[({\bf P}_1,{\bf P}_2)]$ are determined by all Euler characteristics $\chi({\bf P}_1,{\bf P}_2)$ and vice versa, although involving an infinite system of coupled functional equations. To extract more direct information on the relation between these two geometries, we  restrict to more particular cases in the following sections.

\section{Specialization}\label{special}

The system of functional equations in Theorem \ref{gwmk2}, as well as some of the special cases considered in the following sections, simplify considerably once we specialize all variables $s_k$ and $t_l$ to one variable $t$. We denote by $N_{(a,b)}[k]$ the sum
$$N_{(a,b)}[k]=\sum_{|{\bf P}_1|=ka,\, |{\bf P}_2|=kb}N_{(a,b)}[({\bf P}_1,{\bf P}_2)]$$
of Gromov-Witten invariants and by $\chi(k)$ the corresponding sum
$$\chi_{(a,b)}(k)=\sum_{|{\bf P}_1|=ka,\, |{\bf P}_2|=kb}\chi({\bf P}_1,{\bf P}_2)$$
of Euler characteristics. We denote by $f_{(a,b)}(t)\in\Qn [x^ay^b][[t]]$ the specialization of the series $f_{(a,b)}$ (and similarly $R^{{\bf P}_1,{\bf P}_2}(t)$) and define
$$E=\frac{l_1l_2ab-l_2a^2-l_1b^2}{l_1l_2}\in\Qn .$$

\begin{satz}\label{fes} The series $f_{(a,b)}(t)$ is determined by the single functional equation
$$f_{(a,b)}(t)=\prod_{k\geq 1}(1-((tx)^a(ty)^bf_{(a,b)}(t)^E)^k)^{-k\chi_{(a,b)}(k)}.$$
\end{satz}
{\it Proof.} Combining the functional equations of Theorem \ref{gwmk2}, we have $f_{(a,b)}(t)=$
$$\prod_{k\geq 1}\prod_{\substack{{|{\bf P}_1|=ka}\\{|{\bf P}_2|=kb}}}(1-(tx)^{ka}(ty)^{kb}\prod_{k'\geq 1}\prod_{\substack{{|{\bf P}_1'|=k'a}\\{|{\bf P}_2'|=k'b}}}R^{{\bf P}_1',{\bf P}_2'}(t)^{-\langle d({\bf P}_1,{\bf P}_2),d({\bf P}_1',{\bf P}_2')\rangle \chi({\bf P}_1',{\bf P}_2')})^{-k\chi({\bf P}_1,{\bf P}_2)}.$$
We have to study the inner product
$$\prod_{\substack{{|{\bf P}_1'|=k'a}\\{|{\bf P}_2'|=k'b}}}R^{{\bf P}_1',{\bf P}_2'}(t)^{-\langle d({\bf P}_1,{\bf P}_2),d({\bf P}_1',{\bf P}_2')\rangle \chi({\bf P}_1',{\bf P}_2')}$$
on the right hand side.\\[1ex]
Using the $S_{l_1}\times S_{l_2}$-symmetry of the quiver, we see that both the series $R^{{\bf P}_1,{\bf P}_2}(t)$ and the Euler characteristic $\chi({\bf P}_1,{\bf P}_2)$ are invariant under permutations. Thus, in the above product over (pairs of) ordered partitions, it suffices to multiply contributions from ordered partitions ${\bf P}$ satisfying $p_1\geq\ldots\geq p_l$, and to take multiplicities into account. This allows us to simplify the Euler form:\\[1ex]
Let ${\bf P}^0$ be an ordered partition of length $l$ satisfying $p_1^0\geq\ldots\geq p_l^0$, and denote by $z({\bf P}^0)$ the number of rearrangements of ${\bf P}^0$ into an ordered partition ${\bf P}'\models {\bf P}^0$. For an ordered partition ${\bf P}$ of length $l$, we then have

$$\sum_{{\bf P}'\models{\bf P}^0}\sum_kp_kp_k'=\frac{z({\bf P}^0)}{l!}\sum_{\sigma\in S_l}\sum_kp_kp^0_{\sigma(k)}=$$
$$=\frac{z({\bf P}^0)}{l!}\sum_kp_k(\sum_{\sigma\in S_l}p^0_{\sigma(k)})=\frac{z({\bf P}^0)}{l!}(l-1)!\sum_{k,k'}p_kp^0_{k'}=z({\bf P}^0)\frac{1}{l}|{\bf P}||{\bf P}^0|.$$
This implies the following identity for the Euler form:
$$\sum_{\substack{{\bf P}_1'\models{\bf P}_1^0}\\{{\bf P}_2'\models{\bf P}_2^0}}\langle d({\bf P}_1,{\bf P}_2),d({\bf P}_1',{\bf P}_2')\rangle={z({\bf P}_1^0)z({\bf P}_2^0)}kk'\frac{l_2a^2+l_1b^2-l_1l_2ab}{l_1l_2},$$
which is just $-z({\bf P}_1^0)z({\bf P}_2^0)kk'E$.
Using this identity, the inner product above simplifies as follows (${\bf P}_1^0$, ${\bf P}_2^0$ denoting weakly descending ordered partitions as before):
$$\prod_{\substack{{|{\bf P}_1'|=k'a}\\{|{\bf P}_2'|=k'b}}}R^{{\bf P}_1',{\bf P}_2'}(t)^{-\langle d({\bf P}_1,{\bf P}_2),d({\bf P}_1',{\bf P}_2')\rangle \chi({\bf P}_1',{\bf P}_2')}=$$
$$=\prod_{\substack{{|{\bf P}^0_1|=k'a}\\{|{\bf P}^0_2|=k'b}}}\prod_{\substack{{{\bf P}_1'\models{\bf P}^0_1}\\{{\bf P}_2'\models{\bf P}^0_2}}}R^{{\bf P}^0_1,{\bf P}^0_2}(t)^{-\langle d({\bf P}_1,{\bf P}_2),d({\bf P}_1',{\bf P}_2')\rangle \chi({\bf P}_1^0,{\bf P}_2^0)}=$$
$$=\prod_{\substack{{|{\bf P}^0_1|=k'a}\\{|{\bf P}_2^0|=k'b}}}R^{{\bf P}^0_1,{\bf P}^0_2}(t)^{-\sum_{{{\bf P}_1'\models{\bf P}^0_1},\;{{\bf P}_2'\models{\bf P}^0_2}}\langle d({\bf P}_1,{\bf P}_2),d({\bf P}_1',{\bf P}_2')\rangle \chi({\bf P}_1^0,{\bf P}_2^0)}=$$
$$=\prod_{\substack{{|{\bf P}^0_1|=k'a}\\{|{\bf P}_2^0|=k'b}}}R^{{\bf P}^0_1,{\bf P}^0_2}(t)^{{z({\bf P}^0_1)z({\bf P}^0_2)}kk'E \chi({\bf P}_1',{\bf P}_2')}=$$
$$=\prod_{\substack{{|{\bf P}_1'|=k'a}\\{|{\bf P}_2'|=k'b}}}R^{{\bf P}_1',{\bf P}_2'}(t)^{kk'E \chi({\bf P}_1',{\bf P}_2')}.$$

Rewriting the above expression for $f_{(a,b)}(t)$ using this, we get $f_{(a,b)}(t)=$
$$\prod_{k\geq 1}\prod_{\substack{{|{\bf P}_1|=ka}\\{|{\bf P}_2|=kb}}}(1-(tx)^{ka}(ty)^{kb}\prod_{k'\geq 1}\prod_{\substack{{|{\bf P}_1'|=k'a}\\{|{\bf P}_2'|=k'b}}}R^{{\bf P}_1',{\bf P}_2'}(t)^{kk'E\chi({\bf P}_1',{\bf P}_2')})^{-k\chi({\bf P}_1,{\bf P}_2)}.$$
We identify the original product expansion for the series $f_{(a,b)}(t)$ in the inner product, which yields the claimed functional equation.
\qed

We can now give a formula for the specialized Gromov-Witten invariants $N_{(a,b)}[k]$ applying the methods of \cite{ReFE} to the functional equation of Theorem \ref{fes}:

\begin{kor}\label{sgw} We have
$$N_{(a,b)}[k]=\frac{1}{Ek^2}\sum_{{\bf r}}\prod_i\binom{Eki\chi_{(a,b)}(i)+r_i-1}{r_i},$$
the sum running over all ordered partitions ${\bf r}=r_1+\ldots$ such that $\sum_iir_i=k$.
\end{kor}
{\it Proof.} We apply \cite[Proposition 4.4]{ReFE} to the series $F=f_{(a,b)}(t)^E$. Unwinding the definitions, the formula follows.
\qed


\section{Coprime case}\label{coprime}

As a first step towards extracting explicit formulas out of the functional equations of Theorem \ref{gwmk2} (for non-specialized variables), we can compare the coefficients of $x^ay^b$ in the equation of Theorem \ref{gwmk1} as well as in the functional equations of Theorem \ref{gwmk2} to get:

\begin{kor}\label{cor1} For coprime $(a,b)$ and a pair of ordered partitions $({\bf P}_1,{\bf P}_2)$ such that $|{\bf P}_1|=a$ and $|{\bf P}_2|=b$, we have

$$N_{(a,b)}[({\bf P}_1,{\bf P}_2)]=\frac{1}{b}\cdot \chi(M^{\bf b}({\bf P}_1,{\bf P}_2))=\frac{1}{a}\cdot\chi(M^{\bf f}({\bf P}_1,{\bf P}_2))$$
and
$$N_{(a,b)}[({\bf P}_1,{\bf P}_2)]=\chi(M^{\rm st}({\bf P}_1,{\bf P}_2)).$$
\end{kor}

For general coprime $(a,b)$, the Euler characteristic $\chi({\bf P}_1,{\bf P}_2)$ can be computed (preferably with computer aid) by
working out Theorem \ref{rhnr} for the quiver $K(l_1,l_2)$ to obtain the Poincar\'e polynomial, and then specializing $q=1$:

\begin{satz} For a pair $({\bf P}_1,{\bf P}_2)$ of ordered partitions such that $|{\bf P}_1|=a$ and $|{\bf P}_2|=b$ are coprime, we have $\sum_i\dim H^i(M^{\rm st}({\bf P}_1,{\bf P}_2),\Qn )q^{i/2}=(q-1)\times$
$$\times\sum(-1)^{s-1}q^{\sum_{r\leq s}(a_rb_s-\sum_k p_{1,k}^rp_{1,k}^s-\sum_l p_{2,l}^rp_{2,l}^s)}\prod_{r=1}^t(\prod_k\prod_{j=1}^{p_{1,k}^r}(1-q^{-j})\prod_l\prod_{j=1}^{p_{2,l}^r}(1-q^{-j}))^{-1},$$
where the sum runs over all decompositions ${\bf P}_i={\bf P}_i^{(1)}+\ldots+{\bf P}_i^{(t)}$ for $i=1,2$ into ordered partitions ${\bf P}_i^{(r)}=p_{i,1}^r+\ldots+p_{i,l_i}^r$ such that for $a_r=|{\bf P}_1^{(r)}|$ and $b_r=|{\bf P}_2^{(r)}|$, we have $(a_r,b_r)\not=(0,0)$ for all $r=1,\ldots,t$ and
$$\frac{b_1+\ldots+b_r}{a_1+\ldots+a_r}>\frac{b}{a}$$
for all $r<t$. 
\end{satz}

\begin{bem}
\end{bem}
\begin{itemize}
\item In the present context, this formula seems to be somewhat tautological: the factorization formula Theorem \ref{wcqm} ultimately follows from the Harder-Narasimhan recursion for the moduli spaces $M_d^{\Theta-{\rm st}}(Q)$, and the formula in Theorem \ref{rhnr} is a resolution of the very same recursion.\\[1ex]
Also note the disadvantage that the full Poincar\'e polynomial has to be computed to extract just the Euler characteristic, since every individual summand in the above sum has a pole at $q=1$.
\end{itemize}

As a particular example, localization techniques allow to extract the following explicit formula in specialized variables (see Theorem \ref{ddminus1}):

\begin{satz}\label{dm1d} For arbitrary $d$, we have
$$N_{(d,d-1)}[1]=\chi_{(d,d-1)}(1)=\frac{l_1l_2}{d((l_1-1)d+1)}\binom{(l_1-1)(l_2-1)d+l_2-1}{d-1}.$$
\end{satz}

Another particular example is $N_{(3,5)}[1]=204$ for $l_1=3=l_2$, which is worked out in detail in Example \ref{ex184}.

\section{Small length}\label{smalllength}

To work out all series $f_{(a,b)}$ in the case where $l_1l_2\leq 4$, we first note the following trivial case of the functional equations of Theorem \ref{gwmk2}:

\begin{lem} Suppose that the following holds for coprime $a$ and $b$:
\begin{itemize}
\item For all pairs of ordered partitions $({\bf P}_1,{\bf P}_2)$ such that $|{\bf P}_1|=a$ and $|{\bf P}_2|=b$, the moduli space $M^{\rm st}({\bf P}_1,{\bf P}_2)$ is either empty or a single point. Let $(({\bf P}_1^{(r)},{\bf P}_2^{(r)}))_r$ be a complete list of those where the latter holds.
\item The Euler form fulfills $\langle d({\bf P}_1^{(r)},{\bf P}_2^{(r)}),d({\bf P}_1^{(s)},{\bf P}_2^{(s)})\rangle=\delta_{r,s}$.
\item For all pairs of ordered partitions $({\bf P}_1,{\bf P}_2)$ such that $|{\bf P}_1|=ka$ and $|{\bf P}_2|=kb$ for some $k\geq 2$, the moduli space $M^{\rm st}({\bf P}_1,{\bf P}_2)$ is empty.
\end{itemize}
Then
$$f_{(a,b)}=\prod_r(1+s^{{\bf P}_1^{(r)}}t^{{\bf P}_2^{(r)}}x^ay^b),$$
and thus
$$N_{(ka,kb)}[(k{\bf P}_1^{(r)},k{\bf P}_2^{(r)})]=\frac{(-1)^{k-1}}{k^2}.$$
\end{lem}
{\it Proof.} By Theorem \ref{gwmk2}, we have $f_{(a,b)}=\prod_rR^{{\bf P}_1^{(r)},{\bf P}_2^{(r)}}$ and
$$R^{{\bf P}_1^{(r)},{\bf P}_2^{(r)}}=(1-s^{{\bf P}_1^{(r)}}t^{{\bf P}_2^{(r)}}x^ay^b(R^{{\bf P}_1^{(r)},{\bf P}_2^{(r)}})^{-1})^{-1}.$$
The solution to this latter equation is evidently $1+s^{{\bf P}_1^{(r)}}t^{{\bf P}_2^{(r)}}x^ay^b$.
\qed
To complete the description of the series $f_{(a,b)}$ when $l_1l_2\leq 4$, there are thus only two remaining cases to consider:\\[1ex]
We first consider the case $l_1=2=l_2$ and the slope $(a,b)=(1,1)$. We abbreviate $R^{(1+0,1+0)}$ by $R_{11}$ and define $R_{12}$, $R_{21}$ and $R_{22}$ similarly. After computing the relevant values of the Euler form, we find
$$f_{(1,1)}=R_{11}R_{12}R_{21}R_{22},$$
where
$$R_{11}=(1-s_1t_1xy\frac{R_{22}}{R_{11}})^{-1},\;R_{22}=(1-s_2t_2xy\frac{R_{11}}{R_{22}})^{-1},$$
and thus
$$R_{11}=\frac{1+s_1t_1xy}{1-s_1s_2t_1t_2x^2y^2}$$
(and similarly for $R_{12}$, $R_{21}$ and $R_{22}$), thus
$$f_{(1,1)}=\frac{(1+s_1t_1xy)(1+s_1t_2xy)(1+s_2t_1xy)(1+s_2t_2xy)}{(1-s_1s_2t_1t_2x^2y^2)^4}.$$
This gives
$$N_{(1,1)}[(k+0,k+0)]=\frac{(-1)^{k-1}}{k^2}$$
and
$$N_{(1,1)}[(k+k,k+k)]=\frac{2}{k^2}.$$
Similarly, we treat the case $l_1=1$, $l_2=4$ and the slope $(a,b)=(1,2)$: we abbreviate $R^{(1,1+1+0+0)}$ by $R_{12}$ and similarly for $R_{13}$, $R_{14}$, $R_{23}$, $R_{24}$ and $R_{34}$. Then we have
$$f_{(1,2)}=R_{12}R_{13}R_{14}R_{23}R_{24}R_{34}(R^{(2,1+1+1+1)})^{-2}.$$
We have the equation
$$R^{(2,1+1+1+1)}=(1-s^2t_1t_2t_3t_4x^2y^4)^{-1}$$
and
$$R_{12}=(1-st_1t_2xy^2\frac{R_{34}}{R_{12}})^{-1}$$
with solution
$$R_{12}=\frac{1+st_1t_2xy^2}{(1-s^2t_1t_2t_3t_4x^2y^4)^2}$$
as above (and similarly for $R_{13}$, $R_{14}$, $R_{23}$, $R_{24}$ and $R_{34}$). This yields the result
$$f_{(1,2)}=\frac{\prod_{1\leq i<j\leq 4}(1+st_it_jxy^2)}{(1-s^2t_1t_2t_3t_4x^2y^4)^4}.$$
Again, this gives
$$N_{(1,2)}[(k,k+k+0+0)]=\frac{(-1)^{k-1}}{k^2}$$
and
$$N_{(1,2)}[(2k,k+k+k+k)]=\frac{2}{k^2}.$$

\section{Central slope}\label{central}

Now we specialize Theorem \ref{gwmk2} to the central slope $a=1=b$; on the Gromov-Witten side, we are thus considering maps to ${\bf P}^2$ instead of an arbitrary weighted projective plane. We then consider pairs of ordered partitions ${\bf P}_1$, ${\bf P}_2$ such that $|{\bf P}_1|=k=|{\bf P}_2|$. In this case, localization techniques allow us to derive (see Corollary \ref{vanish}) 
that $\chi({\bf P}_1,{\bf P}_2)=0$ as soon as $k\geq 2$. For $k=1$, the choice of ${\bf P}_1$ (resp.~${\bf P}_2$) is just the choice of an index $k=1,\ldots,l_1$ (resp.~$l=1,\ldots,l_2$), and the resulting moduli spaces are single points, thus $\chi({\bf P}_1,{\bf P}_2)=1$. We denote the corresponding dimension vector by $d(k,l)$. The Euler form evaluates to
$$\langle d(k,l),d(k',l')\rangle=\delta_{k,k'}+\delta_{l,l'}-1.$$
Furthermore, we have $$E=\frac{l_1l_2-l_1-l_2}{l_1l_2}$$ in this case. The functional equations of Theorem \ref{gwmk2} then simplify drastically and become algebraic:

\begin{satz}\label{fecs} The series $f_{(1,1)}$ is given as
$$f_{(1,1)}=\prod_{k=1}^{l_1}\prod_{l=1}^{l_2}R^{k,l},$$
where the series $R^{k,l}$ are determined by the system of functional equations
$$R^{k,l}=1+s_kt_lxy\prod_{k'\not=k}\prod_{l'\not=l}R^{k',l'}.$$
\end{satz}

This system of functional equations is reminiscent of the $Q$-systems of \cite{KNT} and can be solved by multivariate Lagrange inversion (but the resulting formulas are not particularly explicit).\\[1ex]
Specializing all variables $s_k$ and $t_l$ to one variable $t$ and defining $$f_{(1,1)}(t)=H(t)^{l_1l_2},$$
the series $H$ is determined by the single functional equation
$$H(t)=(1-t^2xyH(t)^{l_1l_2-l_1-l_2})^{-1}.$$
But then it follows immediately from \cite[Theorem 1.4]{ReHilb} that:

\begin{kor} We have
$$f_{(1,1)}(t)=(\sum_{k\geq 0}\frac{1}{(l_1l_2-l_1-l_2)k+1}\binom{(l_1-1)(l_2-1)k}{k}(t^2xy)^k)^{l_1l_2},$$
\end{kor}

confirming \cite[Conjecture 1.4]{GP}.\\[1ex]
Specializing Corollary \ref{sgw}, we find
\begin{kor} We have
$$N_{(1,1)}[k]=\frac{l_1l_2}{k^2}\binom{(l_1-1)(l_2-1)k-1}{k-1}.$$
\end{kor}
In this special case, we can also confirm (a variant of) the integrality conjecture \cite[Conjecture 6.2]{GPS}:
\begin{kor}  We have
$$\sum_{d|k} \mu(\frac{k}{d})(-1)^{(l_1l_2-l_1-l_2)(d-k)}\frac{d^2}{k^2}N_{(1,1)}[d]\in{\bf N}.$$
\end{kor}
{\it Proof.} By the previous corollary, the above Moebius inversion equals
$$\sum_{d|k} \mu(\frac{k}{d})(-1)^{(l_1l_2-l_1-l_2)(d-k)}\binom{(l_1-1)(l_2-1)d-1}{d-1},$$
which is a nonnegative integer by \cite[Theorem 3.2]{ReQDT}.
\qed

\section{Balanced case}\label{balanced}

In this section, we consider the case $l_1=m=l_2$. We can then relate the geometry of the moduli spaces $M^{\rm st}({\bf P}_1,{\bf P}_2)$ to moduli spaces of representations of the $m$-Kronecker quiver.\\[1ex]
For coprime $a$, $b$ as above and $k\geq 1$, we consider $m$ tuples $(f_k:W\rightarrow V)_k$ of linear maps from a $kb$-dimensional vector space $W$ to a $ka$-dimensional one $V$, up to the base change action of ${\rm GL}(V)\times{\rm GL}(W)$. We call such a tuple of linear maps stable if $\dim\sum_kf_k(U)>\frac{a}{b}\dim U$ for all non-zero proper subspaces $U$ of $W$. There exists a moduli space $M^{\rm st}(ka,kb)$ parametrizing stable tuples up to base change.

From Theorem \ref{kroneckerl1l2}, we derive

\begin{satz}\label{divi} We have
$$\sum_{|{\bf P}_1|=ka,\; |{\bf P}_2|=kb}\chi(M^{\rm st}({\bf P}_1,{\bf P}_2))=m\cdot \chi(M^{\rm st}(ka,kb)).$$
\end{satz}

This allows us to apply \cite[Theorem 5.1]{ReFE} to confirm (again, a variant of) \cite[Conjecture 6.2]{GPS} in the balanced case and for specialized variables, as already indicated in \cite{GP}:

\begin{kor}\label{integ} If $l_1=m=l_2$, every specialized series $f_{(a,b)}(t)$ admits a product factorization
$$f_{(a,b)}(t)=\prod_{k\geq 1}(1-((-1)^{mab-a^2-b^2}t)^k)^{-kd(a,b,k)}$$
for integral $d(a,b,k)$.
\end{kor} 
{\it Proof.} Applying Theorem \ref{fes} and using the above theorem, the series $f_{(a,b)}(t)^{1/m}$ is determined by the functional equation
$$f_{(a,b)}(t)^{1/m}=\prod_{k\geq 1}(1-((tx)^a(ty)^b(f_{(a,b)}^{1/m})^{mab-a^2-b^2})^k)^{-k\chi(M^{\rm st}(ka,kb)}.$$
Applying \cite[Theorem 4.9]{ReFE}, the statement follows.
\qed

\section{General commutator formula}\label{gcf}

In this section, we consider a more general class of bipartite quivers to obtain a partial GW/Kronecker correspondence for the commutator formula \cite[Theorem 5.6]{GPS}, which on the Gromow-Witten side involves orbifold blow-ups. However, it will turn out that the correspondence is weaker than in the cases considered before. The following derivation of the correspondence follows the steps of Section \ref{rc} closely, thus some details will be omitted.\\[2ex]
Fix tuples of nonnegative integers $l_i^*=(l_i^1,\ldots,l_i^{d_i})$ for $i=1,2$. Define a quiver $K(l_1^*,l_2^*)$ as follows: $K(l_1^*,l_2^*)$ has vertices $i_\zeta^r$ for $r=1,\ldots,d_1$ and $\zeta=1,\ldots,l_1^r$ and $j_\xi^s$ for $s=1,\ldots,d_2$ and $\xi=1,\ldots,l_2^s$. There are $rs$ arrows from each $j^s_\xi$ to each $i^r_\zeta$.\\[1ex]
The quiver $K(l_1^*,l_2^*)$ is therefore bipartite with a ``level structure'', the vertices $i^r_\zeta$ and $j^s_\xi$ being of level $r$ and $s$, respectively, such that the number of arrows is given by the product of levels.\\[1ex]
The quiver $K(l_1^*,l_2^*)$ has a natural $\prod_{r=1}^{d_1}S_{l_1^r}\times\prod_{s=1}^{d_2}S_{l_2^s}$-symmetry permuting vertices $i$ (resp.~$j$) of the same level.\\[1ex]
Computing the antisymmetrized Euler form of $K(l_1^*,l_2^*)$, we get
$$\{{i^{r'}_{\zeta'},i^r_\zeta}\}=0,\;\;\; \{ {j^s_\xi,i^r_\zeta}\}=-rs,\;\;\; \{{i^r_\zeta,j^s_\xi}\}=rs,\;\;\; \{{j^{s'}_{\xi'},j^s_\xi}\}=0$$
for all $r,r',s,s',\zeta,\zeta',\xi,\xi'$.\\[1ex]
We view dimension vectors for $K(l_1^*,l_2^*)$ as pairs of graded partitions (in the notation of \cite[Section 5.5]{GPS}) ${\bf G}=({\bf G}_1,{\bf G}_2)$ as follows:\\[1ex]
For $i=1,2$, the graded partition ${\bf G}_i$ is a tuple ${\bf G}_i=({\bf P}^1_i,\ldots,{\bf P}^{d_i}_i)$ of unordered partitions with all parts of ${\bf P}_i^r$ divisible by $r$ (resp.~$s$), that is,
$${\bf P}_1^r=rd_{i_1^r}+\ldots+rd_{i_{l_1^r}^r},\;\;\; {\bf P}_2^s=sd_{j_1^s}+\ldots+sd_{j_{l_2^s}^s}.$$
The collection $(d_{i^r_\zeta},d_{j^s_\xi})$ defines a dimension vector for $K(l_1^*,l_2^*)$. The size of ${\bf G}_i$ is defined as $|{\bf G}_i|=\sum_i|{\bf P}^i|$.\\[2ex]
Specializing the results of Section \ref{qmwcfe} to this quiver, we find the following:\\[1ex]
The algebra $B$ is given by $B=\Qn [[x_{i^r_\zeta},x_{j^s_\xi}]]$; to simplify notation, we rename these variables as $x^{r}_\zeta=x_{i^r_\zeta}$, $y^{s}_\xi=x_{j^s_\xi}$. The Poisson automorphisms $T_{i^r_\zeta}$ and $T_{j^s_\xi}$ of $B$ are given as follows:
$$T_{i^r_\zeta}(x^{r'}_{\zeta'})=x^{r'}_{\zeta'},\;\;\; T_{i^r_\zeta}(y^{s}_\xi)=y^{s}_\xi(1+x^{r}_\zeta)^{rs},$$
$$T_{j^s_\xi}(x^{r}_\zeta)=x^{r}_\zeta(1+y^{s}_\xi)^{-rs},\;\;\; T_{j^s_\xi}(y^{s'}_{\xi'})=y^{s'}_{\xi'}.$$
We consider $B$ as a subalgebra of the algebra $B'=\Qn [x,x^{-1},y,y^{-1}][[s^r_\zeta,t^s_\xi]]$ via the identifications $x^r_\zeta=s^{r}_\zeta x^r$ and $y^s_\xi=t^s_\xi y^s$. Then the automorphisms $T_{i^r_\zeta}$, $T_{j^s_\xi}$ lift to the following $\Qn [[s^r_\zeta,t^s_\xi]]$-linear automorphisms of $B'$:
$$T_{i^r_\zeta}(x)=x,\;\;\; T_{i^r_\zeta}(y)=y(1+s^r_\zeta x^r)^r,$$
$$T_{j^s_\xi}(x)=x(1+t^s_\xi y^s)^{-s},\;\;\; T_{j^s_\xi}(y)=y.$$
We define
$$T_i=\prod_{r=1}^{d_1}\prod_{\zeta=1}^{l^r_1}T_{i^r_\zeta},\;\;\; T_j=\prod_{s=1}^{d_2}\prod_{\xi=1}^{l_2^s}T_{j^s_\xi},$$
thus
$$T_i(x)=x,\;\;\; T_i(y)=y\prod_{r=1}^{d_1}\prod_{\zeta=1}^{l_1^r}(1+s^r_\zeta x^r)^r,$$
$$T_j(x)=x\prod_{s=1}^{d_2}\prod_{\xi=1}^{l_2^s}(1+t^s_\xi y^s)^{-s},\;\;\; T_j(y)=y.$$
We define a $(\Theta,\kappa)$-stability on $K(l_1^*,l_2^*)$ by
$$\Theta_{i^r_\zeta}=0,\;\;\;\Theta_{j^s_\xi}=s,$$
$$\kappa_{i^r_\zeta}=r,\;\;\; \kappa_{j^s_\xi}=s.$$
Given coprime integers $a,b$, we then have the following for a dimension vector ${\bf d}=(d_{i^r_\zeta},d_{j^s_\xi})$ of $K(l_1^*,l_2^*)$:\\[1ex]
$\mu({\bf d})=\frac{b}{a+b}$ if and only if $\sum_r\sum_\zeta rd_{i^r_\zeta}=ka$, $\sum_s\sum_\xi sd_{j^s_\xi}=kb$ for some $k\geq 1$, or, in other words, if and only if $(|{\bf G}_1|,|{\bf G}_2|)$ is a multiple of $(a,b)$.\\[2ex]
By Theorem \ref{wcqm}, we have
$$T_iT_j=\prod^{\leftarrow}_{\mu\in\Qn }T_\mu$$
for automorphisms $T_\mu$ of $B$ given as follows:
$$T_\mu(x^r_\zeta)=x^r_\zeta\prod_{r'}\prod_{\zeta'}(Q_\mu^{i_{\zeta'}^{r'}})^{\{{i^{r'}_{\zeta'},i^r_\zeta}\}}\prod_s\prod_\xi 
(Q_\mu^{j^s_\xi})^{\{{j^s_\xi,i^r_\zeta}\}},$$
$$T_\mu(y^s_\xi)=y^s_\xi\prod_r\prod_\zeta(Q_\mu^{i^r_\zeta})^{\{{i^r_\zeta,j^s_\xi}\}}\prod_{s'}\prod_{\xi'}(Q_\mu^{j^{s'}_{\xi'}})^{\{{j^{s'}_{\xi'},j^s_\xi}\}}.$$
Using the computation of the antisymmetrized Euler form of $K(l_1^*,l_2^*)$ above, this simplifies to
$$T_\mu(x^r_\zeta)=x^r_\zeta\prod_s\prod_\xi 
(Q_\mu^{j^s_\xi})^{-rs},$$
$$T_\mu(y^s_\xi)=y^s_\xi\prod_r\prod_\zeta(Q_\mu^{i^r_\zeta})^{rs}.$$
Using this description, we can see that $T_\mu$ extends to an automorphism of $B'$, namely
$$T_\mu(x)=x\prod_s\prod_\xi 
(Q_\mu^{j^s_\xi})^{-s},$$
$$T_\mu(y)=y\prod_r\prod_\zeta(Q_\mu^{i^r_\zeta})^{r}.$$
The analogue of Lemma \ref{identity} for the more general stability $(\Theta,\kappa)$ yields
$$1=Q_\mu^{\Theta-\mu\kappa}=\prod_r\prod_\zeta(Q_\mu^{i^t_\zeta})^{-\frac{b}{a+b}r}\prod_s\prod_\xi(Q_\mu^{j^s_\xi})^{\frac{a}{a+b}s},$$
and thus
$$1=\prod_r\prod_\zeta(Q_\mu^{i^r_\zeta})^{-br}\prod_s\prod_\xi(Q_\mu^{j^s_\xi})^{as}.$$
By coprimality of $a,b$, we can choose $c,d$ such that $ac+bd=1$. We define
$$F_\mu=\prod_r\prod_\zeta(Q_\mu^{i^r_\zeta})^{cr}\prod_s\prod_\xi(Q_\mu^{j^s_\xi})^{ds}.$$
Then we have
$$F_\mu^a=\prod_r\prod_\zeta(Q_\mu^{i^r_\zeta})^{r},\;\;\; F_\mu^{-b}=\prod_s\prod_\xi(Q_\mu^{j^s_\xi})^{-s},$$
thus
$$T_\mu(x)=xF_\mu^{-b},\;\;\; T_\mu(y)=yF_\mu^a$$
is an element of the tropical vertex $H$.\\[1ex]
The monomials in the variables $x^r_\zeta,y^s_\xi$ appearing in the series $Q_\mu$, $F_\mu$ are the
$$\prod_r\prod_\zeta(x^r_\zeta)^{d_{i^r_\zeta}}\prod_s\prod_\xi(y^s_\xi)^{d_{j^s_\xi}}$$
such that
$$\sum_r\sum_\zeta rd_{i^r_\zeta}=ka,\;\;\; \sum_s\sum_\xi sd_{j^s_\xi}=kb$$
for some $k\geq 1$. After embedding into $B'$, these monomials are identified with
$$\prod_r\prod_\zeta((s^r_\zeta)^{d_{i^r_\zeta}}x^{rd_{i^r_\zeta}})\prod_s\prod_\xi((t^s_\xi)^{d_{j^s_\xi}}y^{sd_{j^s_\xi}})=$$
$$=s^{{\bf G}_1}t^{{\bf G}_2}x^{|{\bf G}_1|}y^{|{\bf G}_2|}=s^{{\bf G}_1}t^{{\bf G}_2}(x^ay^b)^k$$
in the notation of \cite[Section 5.6]{GPS}.\\[1ex]
We consider the moduli space
$$M^{\rm b}({\bf G}_1,{\bf G}_2)=M^{(\Theta,\kappa)}_{d({\bf G}_1,{\bf G}_2),n^{\rm b}}(K(l_1^*,l_2^*))$$
for $d({\bf G}_1,{\bf G}_2)_{i^r_\zeta}=\frac{1}{r}p_{1,\zeta}^r$, $d({\bf G}_1,{\bf G}_2)_{j^s_\xi}=\frac{1}{s}p^s_{2,\xi}$, $n^{\rm b}_{i^r_\zeta}=0$, $n^{\rm b}_{j^s_\xi}=s$ and obtain the following analogue of Theorem \ref{gwmk1}:
\begin{satz}\label{moregeneral} In the tropical vertex $H$, we have a factorization
$$T_{(1,0),\prod_r\prod_\zeta(1+s^r_\zeta x^r)^r}T_{(0,1),\prod_s\prod_\xi(1+t^s_\xi y^s)^s}=\prod_{b/a\mbox{ \footnotesize decreasing}}T_{(a,b),F_{a,b}}$$
where
$$F_{a,b}=\sum_{k\geq 0}(\sum_{\substack{{|{\bf G}_1|=ka,}\\{|{\bf G}_2|=kb}}}\chi(M^{\rm b}({\bf G}_1,{\bf G}_2))s^{{\bf G}_1}t^{{\bf G}_2})(x^ay^b)^k.$$
\end{satz}

The commutator formula \cite[Theorem 5.6]{GPS} involves the product
$$T_{(1,0),\prod_r\prod_\zeta(1+s^r_\zeta x^r)}T_{(0,1),\prod_s\prod_\xi(1+t^s_\xi y^s)}$$
on the left hand side (without the additional powers by $r$ and $s$, respectively) and expresses the factorization on the right hand side in terms of Gromov-Witten invariants $N_{(a,b)}[({\bf G}_1,{\bf G}_2)]$ of orbifold blow-ups of the open surface $X^o_{a,b}$. Thus, the corresponding refinement of the GW/Kronecker correspondence is weaker than the versions in the previous sections. Nevertheless, some partial analogues of the localization techniques used in the previous sections are still available (see Section \ref{vanish}). 

\section{Review of localization theory}\label{ecms}
In the following, let $Q$ be a bipartite quiver with vertices $I\cup J$ and $m(i,j)$ arrows between $j\in J$ and $i\in I$. For a quiver $Q$ we denote by \[N_q:=\{q'\in Q_0\mid\exists\alpha:q\rightarrow q'\vee\alpha:q'\rightarrow q\}\]
the set of neighbours of $q$. Moreover, define $m(i,J):=\sum_{j\in J}m(i,j)$ and $m(I,j)$ analogously.\\[1ex]
For a representation $X$ of the quiver $Q$ we denote by $\underline{\dim} X\in\Nn Q_0$ its dimension vector. Moreover, we choose a level $l:Q_0\rightarrow\Nn^+$ on the set of vertices. Define two linear forms $\Theta, \kappa\in\Hom(\Zn Q_0,\Zn)$ by $\Theta(d)=\sum_{j\in J}l(j)d_j$ and $\kappa(d)=\sum_{q\in Q_0}l(q)d_q$.\\[1ex]
Finally, we define a slope function $\mu:\Nn Q_0\rightarrow\Qn$ by \[\mu(d)=\frac{\Theta(d)}{\kappa(d)}.\]
For a representation $X$ of the quiver $Q$ we define $\mu(X):=\mu(\underline{\dim}X)$.
\begin{defi} A representation $X$ of $Q$ is semistable (resp. stable) if for
all proper subrepresentations $0\neq U \subsetneq X$ the following holds:
\[\mu(U)\leq \mu(X)\text{ (resp. }\mu(U) < \mu(X)).\]
\end{defi}
Fixing a slope function as above, we denote by $R^{\Theta-\rm{sst}}_d(Q)$ the set of semistable points and by $R^{\Theta-\mathrm{st}}_d(Q)$ the set of stable points in the affine variety $R_d(Q)$ of representations of dimension $d\in\Nn Q_0$. Moreover, let $M^{\Theta-\rm{st}}_d(Q)$ (resp. $M^{\Theta-\rm{sst}}_d(Q)$) be the moduli space of stable (resp. semistable) representations. Denote by $\chi$ the Euler characteristic in singular cohomology. Note that if $\kappa=\dim$, we obtain the usual definition of stability.\\[1ex]
In this setup it is easy to check that the stability condition is equivalent to 
\[\sum_{i\in I} l(i)d'_i>\frac{\sum_{i\in I} l(i)d_i}{\sum_{j\in J} l(j)d_j}\sum_{j\in J} l(j)d'_j\]
for all subrepresentations of dimension $d'$.\\[1ex]
Following \cite{sch} we say that a general representation of dimension $d$ satisfies some property if there exists a non-empty open subset $U\subseteq R_d(Q)$ that satisfies this property. By $d'\hookrightarrow d$ we denote if a general representation of dimension $d$ has a subrepresentation of dimension $d'$.\\[1ex]
Fix a quiver $\mathcal{Q}$ and two subquivers $\mathcal{Q}^1$ and $\mathcal{Q}^2$ such that $\mathcal{Q}^1\cup\mathcal{Q}^2=\mathcal{Q}$ and $\mathcal{Q}^1\cap\mathcal{Q}^2=\{q\}$
with $q\in\mathcal{Q}_0$. Then the vertex $q$ is called glueing vertex of $\mathcal{Q}^1$ and $\mathcal{Q}^2$. In the following we denote by $\mathcal{Q}=(\mathcal{Q}^1,\mathcal{Q}^2,q)$ if the quiver $\mathcal{Q}$ is obtained by glueing two quivers $\mathcal{Q}^1$ and $\mathcal{Q}^2$ at the vertex $q$.
\begin{bem}\label{stability}
\end{bem}
\begin{itemize}
\item We consider the Dynkin quiver $A_n=(\{q_1,\ldots,q_n\},\{\alpha_i:q_i\rightarrow q_{i+1}\mid i=1,\ldots,n-1\})$. Let $(Q,d)$ be a tuple consisting of a quiver and a dimension vector and let $j\in Q_0$ be a source of level $l(j)$. We call the tuple 
$((Q,A_{l(j)},j=q_{l(j)}),\hat{d})$ a simple extension of $(Q,d)$ at $j$ if $\hat{d}_q=d_q$ for all $q\in Q_0$ and $\hat{d}_q=d_j$ for all $q\in (A_{l(j)})_0$. We proceed analogously for sinks $i$. We denote the tuple obtained by simple extensions at every vertex by $(Q,Q_0,\hat{d})$.\\[1ex]
Now fix some arbitrary level and let $(Q,d)$ be a tuple as above and consider the simple extension  $(Q,Q_0,\hat{d})$. Obviously, every representation of $(Q,d)$ defines a representation of $(Q,Q_0,\hat{d})$ just by defining the corresponding maps to be the identity. On this simple extension we fix the linear form which takes the value $1$ at every vertex induced by a source $j\in Q_0$ and $j$ itself and the value $0$ at every vertex induced by a sink $i\in Q_0$ and $i$ itself. Now it is easy to verify that a representation of $(Q,d)$ is stable if and only if the corresponding representation is stable with respect to the slope function induced by this linear form. Moreover, we get that a representation of dimension $\hat{d}$ is stable if
\[\sum_{i\in I}\hat{d}'_i>\frac{\sum_{i\in I}\hat{d}_i}{\sum_{j\in J} \hat{d}_j}\sum_{j\in J}\hat{d}'_j\]
for all subrepresentation of dimension $\hat{d'}$. In particular, the whole machinery that is known for moduli spaces of quivers applies in this situation.
\end{itemize}
Now we prove the localization theorem in a slightly more general form than \cite[Corollary 3.15]{wei}.

\begin{satz} We have
\[\chi(M^{\Theta-\rm{st}}_d(Q))=\sum_{\tilde{d}}\chi(M^{\tilde{\Theta}-\rm{st}}_{\tilde{d}}(\tilde{Q})),\]
where $\tilde{d}$ ranges over all equivalence classes being compatible with $d$, and the slope function considered on $\tilde{Q}$ is the one induced by the slope function fixed on $Q$.
\end{satz}
{\it Proof.} The only difference from \cite[Section 3]{wei} is the fact that $d$ is not assumed to be $\Theta$-coprime. Inspection of \cite[Section 3]{wei} shows that this assumption is only required in the proof of \cite[Lemma 3.8]{wei}. It remains to strengthen this lemma by proving the following: if $\tilde{V}$ is a stable representation of $\tilde{Q}$, the induced representation $V$ of $Q$ is also stable. Semistability is proved in \cite[Lemma 3.8]{wei}, and the exclusion of proper non-zero subrepresentations of $V$ of the same slope can be proved as in \cite[Proposition 4.2]{ReLoc}.
\qed

We call a finite subquiver $\mathcal{Q}$ of the universal covering quiver $\tilde{Q}$ of $Q$ localization quiver if there exists a dimension vector $d\in\Nn \mathcal{Q}_0$ such that $M^{\Theta-\rm{st}}_d(\mathcal{Q})\neq\emptyset$. Moreover, fixing such a dimension vector corresponding to a non-empty moduli space we call the tuple $(\mathcal{Q},d)$ localization data.
\begin{bem}\label{uncoloured}
\end{bem}
\begin{itemize}
\item A localization data comes along with a colouring of the arrows $c:\mathcal{Q}_1\rightarrow\ Q_1$ such that arrows which have the same sink or source are coloured differently. Obviously, every such colouring of the arrows gives rise to a localization data. We call a localization data without a fixed embedding uncoloured. Fixing a dimension vector $d\in\Nn Q_0$ we denote by $\mathcal{L}_{d}(Q)$ the set of uncoloured localization data of dimension type $d$, i.e. $\sum_{q\in\mathcal{Q}_0(i)}\tilde{d}_q=d_i$ for $(\mathcal{Q},\tilde{d})\in \mathcal{L}_{d}(Q)$ where $\mathcal{Q}_0(i)$ denotes the set of vertices corresponding to $i$. Moreover, for $(\mathcal{Q},\tilde{d})\in\mathcal{L}_d(Q)$ we denote by $c_Q(\mathcal{Q})$ the set of colourings. Then Theorem \ref{euler} can be stated as
\[\chi(M^{\Theta-\rm{st}}_d(Q))=\sum_{(\mathcal{Q},\tilde{d})\in\mathcal{L}_d(Q)}|c_Q(\mathcal{Q})|\chi(M^{\tilde{\Theta}-\rm{st}}_{\tilde{d}}(\mathcal{Q})).\]
\end{itemize}
Let $(\mathcal{Q},d)$ be an uncoloured localization data such that $\mathcal{Q}=(\mathcal{Q}^1,\mathcal{Q}^2,i)$ for some sink $i$. Let $d^1$ and $d^2$ be the corresponding dimension vectors. Let $X$ be a general stable representation of dimension $d$ and $X_k,\,k=1,2,$ be the corresponding subrepresentations of $\mathcal{Q}^k$. Then there exists a short exact sequence
\[\ses{X_1}{X}{\overline{X}_2}\]
where $\overline{X}_2$ is given by $\ses{S_i^{d_i}}{X_2}{\overline{X}_2}$. Let $\underline{\dim}\overline {X}_2=\bigoplus_{k=1}^l d_k^2$ be the canonical decomposition. Let $q(i)\in Q_0$ be the vertex corresponding to $i$. We split up $(1,\ldots,l)$ into $n(i)=m(q(i),J)-|N_{i}|+1$ (possibly empty) disjoint subsets $S_1,\ldots, S_{n(i)}$ with $\cup S_i=(1,\ldots,l)$. Define $d_{S_t}:=\sum_{k\in S_t}d_k^2$ and $\mathcal{Q}^2_t:=\mathrm{supp}(d_{S_t})$. Then we have the following lemma where we refer to \cite{sch} for a more detailed discussion of canonical decompositions of dimension vectors:
\begin{lem}\label{decomp}
The induced tuple $$(\hat{\mathcal{Q}},\hat{d}):=((\mathcal{Q}^1,\bigcup_{k=1}^{n(i)}\mathcal{Q}^2_k,i),d^1+\sum_{k=1}^{n(i)} d_{S_k}^2)$$ is an uncoloured localization data.
\end{lem}
{\it Proof.} Disregarding colourings, we can view $\hat{Q}$ as a subquiver of $\tilde{Q}$ by definition. A general representation of dimension $\underline{\dim}\overline{X}_2$ decomposes into representations of dimensions $d_1^2,\ldots,d_l^2$, say $X_{d_1^2},\ldots,X_{d_l^2}$. Define $X_{S_t}:=\oplus_{k\in S_t}X_{d^2_k}$. We can understand the representation $\tilde{X}$ which is given by the short exact sequence
\[\ses{X_1}{\tilde{X}}{\bigoplus_{k=1}^{n(i)}X_{S_k}}\]
induced by $\ses{X_1}{X}{\overline{X}_2}$ as a representation of $(\hat{\mathcal{Q}},\hat{d})$. Now every subrepresentation of $\tilde{X}$ naturally induces a subrepresentation of $X$ of the same dimension. Thus, since $X$ is stable, $\tilde{X}$ is also stable.
\qed
\begin{bei}
\end{bei}
\begin{itemize}
\item Consider the generalized Kronecker quiver with dimension vector $(3,5)$ and the localization data given by 
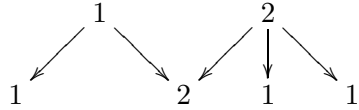
\[
\begin{xy}
\xymatrix@R2pt@C20pt{
&1\ar[lddd]\ar[rddd]&&2\ar[lddd]\ar[ddd]\ar[rddd]&\\\\\\1&&2&1&1}
\end{xy}
\]
Now by applying the lemma to the subquiver on the right hand side we obtain the localization data
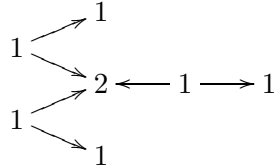
\[\begin{xy}
\xymatrix@R0.5pt@C20pt{
&1\\1\ar[ru]\ar[rd]&\\&2&1\ar[l]\ar[r]&1\\1\ar[ru]\ar[rd]&\\&1 }
\end{xy}
\]
Note that it is not so obvious how to get from the localization data below to the one above.
\end{itemize}
Denote by $S(n)$ the $n$-subspace quiver, i.e. $S(n)_0=\{i,j_1,\ldots,j_n\}$ and $S(n)_1=\{\alpha_k:j_k\rightarrow i\mid k=1,\ldots,n\}$.
\begin{lem}\label{subspace}
Let $d$ be a Schur root of $S(n)$. Then $(S(n),d=(d_0,d_1,\ldots,d_n))$ is an uncoloured localization data of the $n$-Kronecker quiver with dimension vector
$(d_0,\sum_{k=1}^nd_k)$. If we fix a level on the vertices such that $l(i)=\lambda_1$ and $l(j_k)=\lambda_2$ for all $k=1,\ldots,n$ then $(S(n),d=(d_0,d_1,\ldots,d_n))$ is a localization data concerning to the stability induced by this level.
\end{lem}
{\it Proof.}
Following \cite[Theorem 6.1]{sch} for every $d'\hookrightarrow d$ with $d'=(d_0',d'_1,\ldots,d'_n)$ we have
\[\Sc{d'}{d}-\Sc{d}{d'}>0.\]
But it is straightforward to check that this is equivalent to
\[\frac{\sum_{k=1}^{n}d_k}{\sum_{k=0}^nd_k}>\frac{\sum_{k=1}^{n}d'_k}{\sum_{k=0}^nd'_k}\]
which is equivalent to
\[\frac{\lambda_2\sum_{k=1}^{n}d_k}{\lambda_2\sum_{k=1}^nd_k+\lambda_1d_0}>\frac{\lambda_2\sum_{k=1}^{n}d'_k}{\lambda_2\sum_{k=1}^nd'_k+\lambda_1d'_0}.\]
\qed

\section{Localization for $K(l_1,l_2)$}\label{Kl1l2}
We again consider the quiver $K(l_1,l_2)$ where we concentrate on connecting the Euler characteristic of the corresponding moduli spaces to the one of moduli spaces of the Kronecker quiver.\\[1ex] We denote the unique arrow going from $j$ to $i$ by $\alpha_{j,i}$. Without loss of generality we may assume that $l_1\geq l_2$. Let $(\mathbf{P}_1,\mathbf{P}_2)\in\Nn K(l_1,l_2)_0$ and $(a,b)$ its Kronecker type, i.e. $|\mathbf{P}_1|=a$ and $|\mathbf{P}_2|=b$, for not necessarily coprime $a$ and $b$. Define \[\mathcal{L}_{(a,b)}(K(l_1,l_2))=\bigcup_{|{\bf P_1}|=a,|{\bf P_2}|=b}\mathcal{L}_{{\bf P_1},{\bf P_2}}(K(l_1,l_2)).\]
For the Kronecker quiver $K(l_1)$ we denote by $\mathcal{L}^{l_2}_{(a,b)}(K(l_1))$ those uncoloured localization quivers $\mathcal{Q}$ such that $|N_i|\leq l_2$ for all sinks $i\in\mathcal{Q}_0$. Note that the stability conditions on $\tilde{K}(l_1,l_2)$ and $\tilde{K}(l_1)$ coincide. Thus, forgetting the colouring (of the vertices), every $\mathcal{Q}\in \mathcal{L}_{(a,b)}(K(l_1,l_2))$ can be understood as an element of $\mathcal{L}^{l_2}_{(a,b)}(K(l_1))$. Moreover, since there only exists at most one arrow between any two vertices, the set of localization data of dimension type $(\bf{P_1},\bf{P_2})$ of $K(l_1,l_2)$ equals $\mathcal{L}_{{\bf P_1},{\bf P_2}}(K(l_1,l_2))$. Thereby, recall that the vertices of uncoloured localization quivers correspond to vertices of the original quiver. Thus it is straightforward that we have the following lemma:
\begin{lem}\label{1to1}
There exists a one-to-one correspondence between $\mathcal{L}_{(a,b)}(K(l_1,l_2))$ and tuples $(\mathcal{Q},c:\mathcal{Q}\rightarrow K(l_1,l_2)_1)$ where $\mathcal{Q}\in \mathcal{L}^{l_2}_{(a,b)}(K(l_1))$ and $c:\mathcal{Q}_1\rightarrow K(l_1,l_2)_1$ is a colouring such that arrows which have the same sink or source are coloured differently. 
\end{lem}
We get the following statement:
\begin{satz}\label{kroneckerl1l2}
Fix a Kronecker type $(a,b)$ of $K(l_1,l_2)$. Then we have
\[\sum_{|{\bf P}_1|=a,|{\bf P}_2|=b}\chi({\bf P}_1,{\bf P}_2)=\sum_{(\mathcal{Q},\tilde{d})\in\mathcal{L}^{l_2}_{(a,b)}(K(l_1))}|c_{K(l_1,l_2)}(\mathcal{Q})|\chi(M^s_{\tilde{d}}(\mathcal{Q})).\]
If $l_1=l_2=:m$, we have $|c_{K(m,m)}(\mathcal{Q})|=m|c_{K(m)}(\mathcal{Q})|$ for all uncoloured localization data $\mathcal{Q}$. In particular, we have
\[\sum_{{|{\bf P}_1|=a,|{\bf P}_2|=b}}\chi({\bf P}_1,{\bf P}_2)=m\chi(M^s_{a,b}(K(m))).\]
\end{satz}
{\it Proof.} The first statement follows from Lemma \ref{1to1} and the considerations from above. Thus assume that $l_1=l_2$. We choose a map $\mathcal{C}:K(m,m)_1\rightarrow K(m)_1$ such that $\mathcal C(\alpha_{j,i})\neq\mathcal C(\alpha_{j',i})$ and $\mathcal C(\alpha_{j,i})\neq \mathcal C(\alpha_{j,i'})$ for all vertices $i,i',j,j'\in K(m,m)_0$ with $i\neq i'$ and $j\neq j'$.
Every colouring $c:\mathcal{Q}_1\rightarrow K(m,m)_1$ gives rise to a colouring $\mathcal{C}\circ c:\mathcal{Q}_1\rightarrow K(m)_1$. Thus we obtain a surjective map $\mathcal F:c_{K(m,m)}(\mathcal{Q})\rightarrow c_{K(m)}(\mathcal{Q})$. We will proceed by induction on the number of sources of some localization quiver in order to show that $|\mathcal F^{-1}(c)|=m$ for all colourings $c\in c_{K(m)}(\mathcal{Q})$. If $n=1$, the statement is straightforward. Let $\mathcal{Q}$ be a subquiver of $\tilde{K}(m)$ (i.e. we have already fixed some colouring) and $\mathcal{Q}_i,\,i=1,\ldots,m,$ be the corresponding subquivers of $\tilde{K}(m,m)$. If we glue a coloured subquiver of type $(j,i_1,\ldots,i_k)$ with $k\leq m$ to $\mathcal{Q}$, the map $\mathcal{C}$ uniquely determines the colours of the extensions of $\mathcal{Q}_i$. Indeed, the glueing vertex corresponds to one of the $m$ sinks of $K(m,m)$ which only depends on the choice of the previous colouring.
\qed
Note that if $l_2<l_1$ the induction step fails since the map $\mathcal F$ is not surjective in general. Indeed, every sink has only $l_2$ neighbours so that there might be no arrow $\alpha_{j,i}$ that is coloured as needed.\\[1ex]
We consider the quiver $K(l_1,l_2)$ with Kronecker type $(d,d-1)$. Reflecting at every source we may also consider the case $((l_1-1)d+1,d)$ (see also Remark \ref{rem169}). By proceeding as in \cite[Lemma 6.5]{wei}, for all localization data $(\mathcal Q,d)$ we obtain that $d_q=1$ for all $q\in\mathcal Q_0$. Note that if $l_1\geq l_2$ this also follows from the considerations from above. If $l_2>l_1$, the same proof is applicable because it is completely independent of the number of neighbours of some sink. In particular, it only depends on the slope and the number of neighbours of sources.
\begin{satz}\label{ddminus1}
We have
\[\sum_{|{\bf P_1}|=d,|{\bf P_2}|=d-1}\chi({\bf P_1},{\bf P_2})=\frac{l_1l_2}{d((l_1-1)d+1)}\binom{(l_2-1)(l_1-1)d+l_2-1}{d-1}.\]

\end{satz}
{\it Proof.}
We proceed analogously to \cite[Theorem 6.6]{wei}. From the considerations from above we obtain that all sub-localization data of a localization data, which have one source, have vertex set $\{j,i_1,\ldots,i_{l_1}\}$ with $d_j=d_{i_k}=1$. In particular, the moduli spaces of the all considered quivers are points.\\[1ex]
There exists exactly one possibility to colour the arrows of such a quiver taking into account the symmetries of $S_{l_1}$. Now we can glue $k$ subquivers on each vertex $i_l$, $1\leq l\leq l_1$, with $0\leq k\leq (l_2-1)$. But we have to take note of the symmetries of $S_k$. Assuming that there is only one starting vertex to which we can glue, let $y(x)$ be the generating function of such quivers and consider
\begin{eqnarray*}\phi(x)&=&1+\frac{(l_2-1)}{|S_1|}x^{l_1-1}+\frac{(l_2-1)(l_2-2)}{|S_2|}x^{2(l_1-1)}+\ldots+\frac{\prod_{i=1}^{l_2-1}(l_2-i)}{|S_{l_2-1}|}x^{(l_2-1)(l_1-1)}
\\&=&\sum_{i=0}^{l_2-1}x^{i(l_1-1)}\binom{l_2-1}{i}=(1+x^{l_1-1})^{l_2-1}.
\end{eqnarray*}
The generating function of such trees satisfies the functional equation $y(x)=x(\phi(y(x)))$. The generating function for all localization data is obtained as follows: we start with a localization data of Kronecker type $(1,l_1)$ having $l_1$ vertices to which we can glue. The resulting generating function is $y(x)^{l_1}$. By applying the Lagrange inversion theorem, see for instance \cite{sta} for more details, we obtain that
\[[x^n]y(x)^{l_1}=\frac{l_1}{n}[u^{n-l_1}]\phi(u)^n=\frac{l_1}{n}\binom{n(l_2-1)}{\frac{n-l_1}{l_1-1}}.\]
If we assign the weight $0$ to the source of the localization data started with, every such quiver that has $(l_1-1)d+1$ knots corresponds to a localization data of Kronecker type $(d,(l_1-1)d+1)$. The other way around, we may assume that every localization data has some source $j\in J$ with weight $0$ what gives us $d$ choices. This means for every localization data we exactly get $d$ trees. Moreover, we have to take into account that we have $l_2$ choices for the colour of the source of the starting quiver. Hence we get
\[\sum_{|{\bf P_1}|=(l_1-1)d+1,|{\bf P_2}|=d}\chi({\bf P_1},{\bf P_2})=\frac{l_1l_2}{d((l_1-1)d+1)}\binom{(l_2-1)(l_1-1)d+l_2-1}{d-1}.\]
\qed

\section{Vanishing of the Euler characteristic}\label{vanisheuler}
Assume that 
\[\sum_{i\in I}l(i)d_i=K\sum_{j\in J}l(j)d_j\]
for some $K\in\Qn$.
Let $S(n)^t$ be the quiver obtained from the subspace quiver by reversing all arrows.
\begin{defi}
Let $d\in\Nn Q_0$ and $\mathrm{supp}(d)_0=\{i_1,\ldots,i_n,j_1,\ldots,j_m\}$. We say that $d$
satisfies the no-peak condition if the following conditions hold: 
If $m=1$, then there exists no decomposition of the form $$d=(d_j,d_{i_1}^1+\ldots+d_{i_1}^{m(i_1,j)},\ldots,d_{i_n}^1+\ldots+d_{i_n}^{m(i_n,j)})$$ such that
$d=(d_j,d_{i_1}^1,\ldots,d_{i_1}^{m(i_1,j)},\ldots,d_{i_n}^1,\ldots,d_{i_n}^{m(i_n,j)})$ is a Schur root of $S(m(I,j))$.\\
If $n=1$, then there exists no decomposition of the form $$d=(d_i,d_{j_1}^1+\ldots+d_{j_1}^{m(i,j_1)},\ldots,d_{j_m}^1+\ldots+d_{j_m}^{m(i,j_m)})$$ such that
$d=(d_i,d_{j_1}^1,\ldots,d_{j_1}^{m(i,j_1)},\ldots,d_{j_m}^1,\ldots,d_{j_m}^{m(i,j_m)})$ is a Schur root of $S(m(i,J))^t$.
\end{defi}
In other words, if $d$ satisfies the no-peak condition, there does not exist a localization data with only one sink or only one source. For instance if $Q=K(m)$ and $d=(a,a)$, then $d$ satisfies the no-peak condition for all $a\geq 2$.
\begin{pro}\label{neighbours}
For every sink $i$ of some localization quiver with $|I|\geq 2$ we have
\[l(i)<K\sum_{j\in N_i}l(j).\]
For every source $j$ of some localization quiver with $|J|\geq 2$ we have 
\[l(j)<\frac{1}{K}\sum_{i\in N_j}l(i).\]
\end{pro}
{\it Proof.}
Assume that there exists a subquiver 
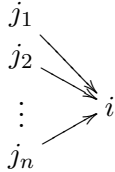
\[
\begin{xy}
\xymatrix@R0.5pt@C20pt{
j_{1}\ar[rdd]&\\j_{2}\ar[rd]&\\\vdots&i\\j_{n}\ar[ru]&}
\end{xy}
\]
By use of
\[\Sc{\underline{\dim} X}{\underline{\dim} Y}=\dim\Hom(X,Y)-\dim\Ext(X,Y)\]
for two representations of a quiver $Q$, see for instance \cite{rin}, we get that a representation of this quiver has a factor of type
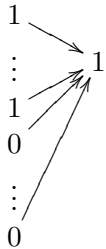
\[
\begin{xy}
\xymatrix@R0.5pt@C20pt{
1\ar[rd]&\\\vdots&1\\1\ar[ru]\\0\ar[ruu]&\\\vdots\\0\ar[ruuuu]}
\end{xy}
\]
if
\begin{eqnarray}\label{1}\sum_{k=1}^sd_{j_k}+d_i-\sum_{k=1}^nd_{j_k}=-\sum_{k=s+1}^nd_{j_k}+d_i>0.\end{eqnarray}
Here $s$ denotes the number of vertices of dimension one on the left hand side of the quiver. Note that we may without lose of generality assume that every homomorphism is surjective because otherwise we would get a surjection for some $t<s$.\\[1ex]
Since such a factor has to be of bigger slope we get
\[l(i)<K\sum_{k=1}^s l(j_k).\]
Since inequality (\ref{1}) definitely holds for $s=n$, the claim follows.\\[1ex]The second claim follows when considering subquivers of type
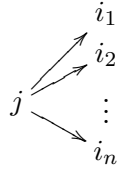
\[
\begin{xy}
\xymatrix@R0.5pt@C20pt{
&i_{1}&\\&i_{2}&\\j\ar[rd]\ar[ruu]\ar[ru]&\vdots\\&i_{n}&}
\end{xy}
\]
and subrepresentations of type
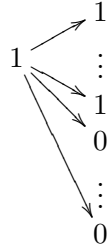
\[
\begin{xy}
\xymatrix@R0.5pt@C20pt{
&1&\\1\ar[ru]\ar[rd]\ar[rdd]\ar[rdddd]&\vdots\\&1\\&0&\\&\vdots\\&0}
\end{xy}
\]
Since we now deal with subrepresentations which must be of smaller slope we obtain
\[Kl(j)<\sum_{k=1}^s l(i_k).\]
\qed
We obtain the following result:
\begin{satz}
Assume that $l(i)=\lambda_1$ for all $i\in I$, $l(j)=\lambda_2$ for all $j\in J$ and let $K=\frac{\lambda_1}{\lambda_2}$. Moreover, assume that
\[\sum_{i\in I}l(i)d_i=K\sum_{j\in J}l(j)d_j\]
and that $d$ satisfies the no-peak condition. Then there does not exist any localization quiver.
\end{satz}
{\it Proof.} 
According to Proposition \ref{neighbours} for every vertex of a localization quiver we have
\[l(i)=\lambda_1<\frac{\lambda_1}{\lambda_2}\sum_{j\in N_i}l(j)=\lambda_1|N_i|\]
and 
\[l(j)=\lambda_2<\frac{\lambda_2}{\lambda_1}\sum_{i\in N_j}l(i)=\lambda_2|N_j|.\]
Thus every vertex of a localization quiver is forced to have at least two neighbours. In particular, every localization quiver is forced to be cyclic. 
\qed
If $l(q)=1$ for all $q\in Q_0$ and $K=1$, this proves a more general version of \cite[Corollary 6.3]{wei}:
\begin{kor}\label{vanish}
Let $d\in\Nn Q_0$ be a dimension vector such that
\[\sum_{i\in I}d_i=\sum_{j\in J}d_j,\]
$l(q)=1$ for all $q\in Q_0$ and $\sum_{j\in J}d_j\neq 1$. Then we have $\chi(M_d^{\Theta-\rm{st}}(Q))=0$.
\end{kor}
Moreover, we get the following result:
\begin{satz}\label{eulervan}
Let $d\in\Nn Q_0$. Assume that
\[\sum_{i\in I}l(i)d_i=K\sum_{j\in J}l(j)d_j\]
and that $d$ satisfies the no-peak condition. Then there exists no localization quiver if 
\[l(i)\geq K\sum_{j\in N_i}m(i,j)l(j)\]
for some sink $i\in I$ or
\[l(j)\geq \frac{1}{K}\sum_{i\in N_j}m(i,j)l(i)\]
for some source $j\in N_i$. In particular, the Euler characteristic of the corresponding moduli space vanishes.
\end{satz}
{\it Proof.}
The claim follows because every neighbour of some vertex of the original quiver gives rise to $m(i,j)$ neighbours in the universal cover.
\qed
If $Q$ has the generalized Kronecker quiver $K(m)$ as a proper subquiver such that the corresponding sink $i$ has only one neighbour, we immediately get the following corollary:
\begin{kor}
Let $d\in\Nn Q_0$ satisfy the no-peak condition. Assume that $K=1$ and that there exist $i,j\in Q_0$ such that $m(i,j)=m$ and $N_i=\{j\}$. If $l(j)=1$ and $l(i)=m$, there does not exist any localization quivers.
\end{kor}
The last case we treat is the following:
\begin{satz}
Let $d\in\Nn Q_0$ satisfy the no-peak condition and assume that $l(i)=l$ for all $i\in I$ and
\[l\sum_{i\in I}d_i=K\sum_{j\in J}l(j)d_j\]
where $K\in\Nn$. Moreover, for all $j\in J$ and $i\in I$ let
\[Kl(j)\geq l\left(m(I,j)-1\right).\]
Then there exists no localization quiver.
\end{satz}
{\it Proof.} Because of Theorem \ref{eulervan} we can assume that
\[Kl(j)= l\left(m(I,j)-1\right)\]
for all vertices $i\in I$ and $j\in J$.\\[1ex]Assume that there exists a localization quiver and let $X$ be a stable representation. Consider a subrepresentation  
\[
\begin{xy}
\xymatrix@R5pt@C20pt{
&X_{i_{1}}&\\&X_{i_{2}}&\\X_j\ar[rd]_{X_{N}}\ar[ruu]^{X_{1}}\ar[ru]_{X_{2}}&\vdots\\&X_{i_N}&}
\end{xy}
\]
such that $|N(i_k)|=1$ for all $k\geq 2$. Then we obviously have $d_{i_k}\leq d_j$ for all $k\geq 2$. Consider the kernel of the map $X_k$. We have
\begin{eqnarray*}\dim(\sum_{t=1}^NX_t(\ker(X_k)))&\leq& \dim(\ker(X_k))(N-1)l\leq\\
&\leq&\dim(\ker(X_k))Kl(j).
\end{eqnarray*}
Thus we obtain $\ker(X_k)=\{0\}$ and, therefore, $d_j=d_{i_k}$ for all $k=2,\ldots,N$. Thus by Lemma \ref{decomp} we may assume that $d_j=d_{i_k}=1$.\\[1ex]
But we also have that $N=m(I,j)$. Indeed, otherwise we would have 
\[\dim(\sum_{l=1}^NX_l(X_j))=Nl\leq l\left(m(I,j)-1\right)=Kl(j).\]
But this means that the representation has a factor representation of dimension type $(l(j),Kl(j))$ which contradicts the stability condition. 

\qed
Moreover, we get the following corollary:
\begin{kor}
Assume that $l(q)=1$ for all vertices of the quiver $Q$ and that
\[K\sum_{j\in J}d_j=\sum_{i\in I}d_i\]
where $\sum_{j\in J}d_j\neq 1$ and $K\in \Nn$.
If we have $m(I,j)\leq K+1$ for all $j\in J$, there exists no localization quiver. In particular, the Euler characteristic of the corresponding moduli space vanishes.
\end{kor}

\begin{bem}\label{rem169}
\end{bem}
\begin{itemize}
\item As far as the Kronecker quiver $K(m)$ is concerned the preceding statement says that the Euler characteristic vanishes if $(d,e)=(d,kd)$ and $k\geq m-1$ which also follows from Theorem \ref{vanish} by applying the reflection functor, see \cite{bgp} for the definition. In the case of bipartite quivers such that $l(q)=1$ for all $q$, the reflection functor applied to all sinks or sources simultaneously gives rise to isomorphisms between moduli spaces. Indeed, it is checked easily that the stability conditions are equivalent. Moreover, subrepresentations become factor representations and vice versa.\\[1ex]
But in the general case, it is not obvious how to get isomorphisms between moduli spaces corresponding to different dimension vectors (except the one coming from transposing all maps) because the stability conditions are not compatible.  
\end{itemize}
\section{Further examples}\label{appl}
In this section we give several examples and applications illustrating the results of the preceding sections. 
\subsection{The case K=1}
Assume that we have
\[\sum_{i\in I} l(i)d_i=\sum_{j\in J}l(j)d_j.\]
Let $(\mathcal{Q},d)$ be a localization data and let
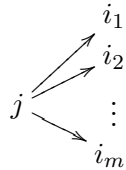
\[
\begin{xy}
\xymatrix@R0.5pt@C20pt{
&i_{1}\\&i_{2}\\j\ar[ruu]\ar[ru]\ar[rd]&\vdots\\&i_{m}\;}
\end{xy}
\]
be a subquiver such that $N_{i_k}=\{j\}$ for all $k=2,\ldots,m$. Because of the stability condition we have
\[\sum_{k=2}^ml(i_k)d_{i_k}<l(j)d_j.\]
In particular, if $l(j)=1$, by Lemma \ref{decomp} it follows that $d_{i_k}=0$ for all $k=2,\ldots,m$. Indeed, in the canonical decomposition of $(d_j,d_{i_2},\ldots,d_{i_k})$ the simple representation $S_j$ at least occurs with multiplicity $r:=d_j-\sum_{k=2}^m l(i_k)d_{i_k}$. In particular, we also get a factor of slope $(d_j-r,\sum_{k=2}^m l(i_k)d_{i_k})=(d_j-r,d_j-r)$ which contradicts the stability condition.\\[1ex]
Recall also Proposition \ref{neighbours} which for $K=1$ says how many neighbours of which kind are allowed for a vertex with a fixed level.\\[1ex]
Now assume that there exists a subquiver
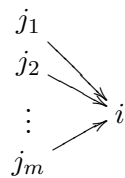
\[
\begin{xy}
\xymatrix@R0.5pt@C20pt{
j_{1}\ar[rdd]&\\j_{2}\ar[rd]&\\\vdots&i\\j_{m}\ar[ru]&}
\end{xy}
\]
such that $l(i)=1=l(j_k)$ for all $k$ and $N_i=\{j_1,\ldots,j_m\}$. Then we have
\[(\sum_{k=1}^md_{j_k})-d_{j_l}\geq d_i\]
for all $l=1,\ldots,m$ because otherwise a representation of this dimension would have a factor isomorphic to $1\rightarrow 1$. Moreover, it again follows that we have $|N_{j_k}|\geq 2$ for all $k$.\\[1ex]
In the following we assume that every vertex has enough neighbours in order to show that in general it can not be assumed that there exist no localization data.\\[1ex]
We restrict to the cases $l(j)=1$ for all $j\in J$ and $l(i)\neq 1$ for all $i\in I$ and assume that $\sum_{j\in J}d_j=\sum_{i\in I}l(i)d_i$ with $d_i\geq 2$ for all $i\in I$. First we deal with the Kronecker quiver with dimension vector $(e,d)$ and $l(i)e=d$. 
Then the stability condition is 
\[l(i)e'>\frac{l(i)e}{d}d'\Leftrightarrow e'>\frac{e}{d}d'\]
for all $(d',e')\hookrightarrow (d,e)$. Thus we can consider the localization quiver
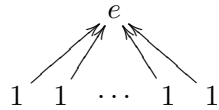
\[
\begin{xy}
\xymatrix@R20pt@C5pt{&&e&&\\1\ar[rru]&1\ar[ru]&\ldots&1\ar[lu]&1\ar[llu]}
\end{xy}
\]
with $d$ vertices on the bottom row. Since we have $d>e$, it is easy to check that this is a Schur root of the subspace quiver. Note that if $l(i)=1$ there only exist semistable and polystable points respectively. But nevertheless, there can exist localization quivers in the case of bipartite quivers with more than one sink and such that $l(i)=1$ not for all sinks of the quiver, see also Example \ref{bei1}.\\[1ex]
Now consider two dimension vectors $(d_i,e_i)$ satisfying the conditions from above. Then we can glue the mentioned localization quivers in order to get one of type $(d_1+d_2,e_1+e_2)$ as follows
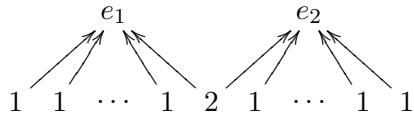
\[
\begin{xy}
\xymatrix@R20pt@C5pt{&&e_1&&&&e_2\\1\ar[rru]&1\ar[ru]&\ldots&1\ar[lu]&2\ar[llu]\ar[rru]&1\ar[ru]&\ldots&1\ar[lu]&1\ar[llu]}
\end{xy}
\]
Now it is easy to check that a general representation of this quiver and with this dimension vector is stable. Moreover, we can recursively construct such quivers. Thus we get localization data for all quivers as above such that $\sum_{j\in J}d_j=\sum_{i\in I}l(i)d_i$ and $l(i)\neq 1$ for all $i\in I$.
\begin{bei}\label{bei1}
\end{bei}
Consider the quiver having vertices $Q_0=\{j,i_1,i_2\}$ with $l(j)=1=l(i_1)$ and $l(i_2)=l$. Moreover, consider the dimension vector $(4l+2,2,4)$ and the localization data
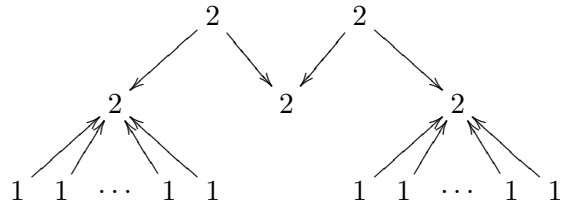
\[
\begin{xy}
\xymatrix@R20pt@C5pt{&&&&2\ar[lld]\ar[rrd]&&&&2\ar[lld]\ar[rrd]\\&&2&&&&2&&&&2\\1\ar[rru]&1\ar[ru]&\ldots&1\ar[lu]&1\ar[llu]&&&&1\ar[rru]&1\ar[ru]&\ldots&1\ar[lu]&1\ar[llu]}
\end{xy}
\]
with twice $2l-1$ vertices on the bottom row and where the single vertex in the middle row has level one. Now it is easily verified that a general representation of this dimension is stable.
\begin{bei}
\end{bei}
The vanishing of the Euler characteristic does not imply that the moduli space of the original quiver is empty. For instance let $Q$ be the quiver with vertices $J=\{j\}$, $I=\{i_1,i_2\}$ and $m(i_1,j)=3$ and $m(i_2,j)\geq 3$. Consider the dimension vector $d$ given by $d_j=21,\,d_{i_1}=6$ and $d_{i_2}=3$ and the level given by $l(j)=1,\,l(i_1)=3$ and $l(i_2)=1$. Thus we have $l(i_1)=m(i_1,j)l(j)$ and by Theorem \ref{eulervan} it follows that $\chi(M_d^{\Theta-\rm{st}}(Q))=0$. But a general representation of the following subquiver of the universal abelian (!) covering quiver
\[
\begin{xy}
\xymatrix@R25pt@C0.5pt{&&1\ar[lld]^{\beta_1}\ar[rrd]^{\alpha_1}&&&&1\ar[lld]^{\alpha_2}\ar[rrd]^{\alpha_3}&&&&1\ar[lld]^{\alpha_1}\ar[rrd]^{\beta_2}&&&&1\ar[lld]^{\beta_3}\ar[rrd]^{\alpha_3}&&&&1\ar[lld]^{\alpha_1}\ar[rrd]^{\alpha_2}&&&&1\ar[lld]^{\alpha_3}\ar[rrd]^{\beta_1}&&&&1\ar[lld]^{\beta_2}\ar[rrd]^{\alpha_2}&&&&1\ar[lld]^{\alpha_3}\ar[rrd]^{\alpha_1}\\1&&&&\underline{1}&&&&\underline{1}&&&&1&&&&\underline{1}&&&&\underline{1}&&&&1&&&&\underline{1}&&&&\underline{1}\\&&&1\ar[ru]&&1\ar[lu]&&1\ar[ru]&&1\ar[lu]&&&&&&1\ar[ru]&&1\ar[lu]&&1\ar[ru]&&1\ar[lu]&&&&&&1\ar[ru]&&1\ar[lu]&&1\ar[ru]&&1\ar[lu]\\\\&&&&&&&&&&&&&&&&1\ar[rrrrrrrrrrrrrrrruuu]^{\alpha_2}\ar[lllllllllllllllluuu]^{\beta_3}&&&&}
\end{xy}
\]
is stable. Here we can choose every appropriate colouring such that different vertices have different weights. Moreover, we denote the arrows between $j$ and $i_1$ by $\alpha_k$, the arrows between $j$ and $i_2$ by $\beta_k$ and by $\underline{~}$ we denote the vertices corresponding to $i_1$. Now it is easy to verify that a general representation of this localization data (corresponding to the universal abelian cover) is stable. \\[1ex]Note that the number of arrows $m(j,i_1)$ plays an important role. If we choose $m(j,i_1)=5$ in the same setup, there exists for instance the following localization data
\[
\begin{xy}
\xymatrix@R20pt@C5pt{&&\underline{3}&&&3&&&\underline{3}&&\\1\ar[rru]&2\ar[ru]&2\ar[u]&2\ar[lu]&3\ar[llu]\ar[ru]&&3\ar[lu]\ar[rru]&2\ar[ru]&2\ar[u]&2\ar[lu]&2\ar[llu]&}
\end{xy}
\]
Note that, the dimension vector $(3,1,2,2,2,2,2,2,2)$ is a Schur root of the $8$-subspace quiver. Thus, it is straightforward to check that a general representation of this dimension is stable.

\subsection{The case $Q=K(l_1,l_2)$}\label{thecase}
We consider the case $m(i,j)=1$ and $l(q)=1$ for all $i,j,q\in Q_0$. Let $d\in\Nn Q_0$ and define $|J|=l_2$, $|I|=l_1$ and $b=\sum_{j\in J}d_j$, $a=\sum_{i\in I}d_i$. 
We may assume that $l_1\geq l_2$. Otherwise we turn around all arrows. Moreover, we can assume that $\frac{l_1}{2}b\geq a$. Otherwise we apply the reflection functor to every sink and turn around all arrows afterwards. Note that, for the reflected Kronecker type $(b,l_1b-a)$ we have
\[\frac{l_1b-a}{b}<\frac{a}{b}\Leftrightarrow \frac{l_1}{2}b<a.\]
Let $d$ be a root of $Q$ and let $m:=\lceil\frac{l_2+l_1}{2}\rceil$.
By an easy calculation we get that $f(d_1,\ldots,d_{l_2}):=\sum_{j\in J}d_j^2$ has a unique minimum at $d_j=\frac{b}{l_2}$, $j\in J$. Thus we get
\begin{eqnarray*}\Sc{d}{d}&=&\sum_{j\in J}d_j^2+\sum_{i\in I}d_i^2-ab\geq\frac{b^2}{l_2}+\frac{a^2}{l_1}-ab=\frac{1}{l_1l_2}(l_1b^2+l_2a^2-l_2l_1ab)\\
&=&\frac{1}{l_2}(b^2+\frac{l_2}{l_1}a^2-l_2ab)\geq \frac{2}{l_2+l_1}(b^2+a^2-\frac{l_2+l_1}{2}ab)
\end{eqnarray*}
where the last inequality holds because $l_1\geq l_2$ and 
\[\frac{l_2}{l_1}a^2-l_2ab\geq a^2-\frac{l_2+l_1}{2}ab\Leftrightarrow\frac{l_1}{2}b\geq a.\]
If $\Sc{d}{d}\leq 0$, it obviously follows that $(b,a)$ is a root of $K(m)$.\\[1ex]
Now we consider the case $l_1=l_2=:m$ and $\Sc{d}{d}=1$ in greater detail. We restrict to the case $m\neq 1$. Then we have $m\nmid b$ or $m\nmid a$ because otherwise we would have
\[(k^2_1m^3+k_2^2m^3-m^4k_1k_2)=m^3(k_1^2+k_2^2-mk_1k_2)=m^2\Sc{d}{d}\] 
for some $k_1,\,k_2\in\Nn$. Thus we get that the first inequality is proper and we may assume that $m\nmid b$. Then we have $d_j=\frac{b+s_j}{m}$ with $s_j\neq 0$ and $\sum_{j=1}^m s_j=0$ because $\sum_{j=1}^md_j=\sum_{j=1}^md_j+s_j=b$. We obtain
\[\sum_{j=1}^md_j^2=\sum_{j=1}^m\left(\frac{b+s_j}{m}\right)^2
=\frac{b^2}{m}+\sum_{j=1}^m\frac{2bs_j+s_j^2}{m^2}=\frac{b^2}{m}+\sum_{j=1}^m\frac{s_j^2}{m^2}.\]
Now we can assume that $s_j>0$ for $j=1,\ldots,n$ and $s_j<0$ for $j=n+1,\ldots,m$ with $n\geq 1$. Define $n$ by $\frac{n}{m}=\lceil\frac{b}{m}\rceil-\frac{b}{m}$. Then we have
$\frac{m-n}{m}=\frac{b}{m}-\lfloor\frac{b}{m}\rfloor$. If we choose $s_j=n-m$ for all $j=1,\ldots, n$ and $s_j=n$ for all $j=n+1,\ldots,m$, we get the unique tuple satisfying $\sum_{j=1}^m s_j=0$ and $|s_j|\leq m$ for every $j$. Because 
\[(s_k-m)^2+(s_l+m)^2=s_k^2+s_l^2-2s_km+2s_lm+2m^2\]
and
\[s_lm-s_km+m^2\geq 0\Leftrightarrow s_k+m\geq s_l,\]
it follows that $\sum_{j=1}^ms_j^2$ takes its minimum at this tuple.\\[1ex]Keeping in mind that $\sum_{j=1}^ns_j=-\sum_{j=n+1}^ms_i$ as before we get
\begin{eqnarray*}\sum_{j=1}^ns_j^2+\sum_{j=n+1}^ms_j^2&\geq &n\left(\frac{\sum_{j=1}^ns_j}{n}\right)^2+(m-n)\left(\frac{\sum_{j=n+1}^ms_j}{m-n}\right)^2\\&=&\frac{m\left(\sum_{j=1}^ns_j\right)\left(-\sum_{j=n+1}^ms_j\right)}{n(m-n)}
\geq\frac{mn^2(m-n)^2}{n(m-n)}\geq m(m-1)
\end{eqnarray*}
where the last inequality holds because $m\geq n+1$.
Thus we get
\[1=\Sc{d}{d}=\frac{b^2}{m}+\sum_{i=1}^m\frac{2bs_j+s_j^2}{m^2}+\frac{a^2}{m}+\sum_{j=1}^m\frac{2at_j+t_j^2}{m^2}-ab\geq \frac{m(m-1)}{m^2}+\frac{\Sc{(b,a)}{(b,a)}}{m}\]
and, therefore,
\[\Sc{(b,a)}{(b,a)}\leq m-(m-1)=1.\]
Thus $(b,a)$ is a root of $K(m)$. We get the following statement:
\begin{lem}
If $({\bf P_1},{\bf P_2})$ is a root of Kronecker type $(b,a)$, then $(b,a)$ is a root of $K(m)$. If $M^{\Theta-\rm{st}}_{(b,a)}(K(m))^T\neq\emptyset$, then there exists a dimension vector $({\bf P'_1},{\bf P'_2})$ of Kronecker type $(b,a)$ such that $M_{({\bf P'_1},{\bf P'_2})}^{\Theta-\rm{st}}(K(m,m))\neq\emptyset$.
\end{lem}
{\it Proof.} The first statement follows from the considerations from above. Thus let $(b,a)$ be a root. By \cite[Theorem 3.11]{wei} torus fixed point sets are described as moduli spaces of the universal abelian covering quiver. Now applying methods similar to the one of Section \ref{Kl1l2} we can construct stable representations of $K(l_1,l_2)$ as claimed.
\qed
Note that even for the universal abelian cover there do not exist stable representations for every dimension vector of $K(m)$, e.g. if $m=3$ and $(d,e)=(2,2)$ or $(d,e)=(4,4)$.
\begin{bei}\label{ex184}
\end{bei}
Consider the example $(3,5)$ of $K(3)$ with the following localization data: 
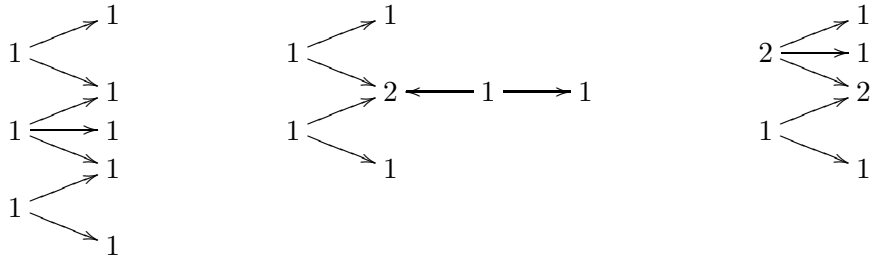
\[
\begin{xy}
\xymatrix@R1.5pt@C25pt{
&1&&&1&&&&&1\\1\ar[ru]\ar[rd]&&&1\ar[ru]\ar[rd]&&&&&2\ar[ru]\ar[rd]\ar[r]&1\\&1&&&2&1\ar[l]\ar[r]&1&&&2\\1\ar[r]\ar[ru]\ar[rd]&1&&1\ar[ru]\ar[rd]&&&&&1\ar[ru]\ar[rd]\\&1&&&1&&&&&1 \\1\ar[ru]\ar[rd]&\\&1}
\end{xy}
\]
In summary, we have $68$ possible colourings of the arrows (i.e. possibilities to embed these quivers into $\tilde{K}(3)$) and, therefore, we get $\chi(M_{(3,5)}^s(K(3)))=68$.\\[1ex]Now consider $K(3,3)$ and ${\bf P_1}=3+1+1$ and ${\bf P_2}=1+1+1$. In order to obtain some localization data we can think of colouring the vertices of the above quivers such that vertices having the same sink or source do not have the same colour. Obviously, this is just possible for the first localization quiver where we have $12$ possibilities (six on the left hand side and two on the right hand side). If we mod out the action of $S_2$ we get $\chi(3+1+1,1+1+1)=6$.\\[1ex]
Now consider ${\bf P_1}=2+2+1$ and ${\bf P_2}=1+1+1$. It is straightforward that we have 18 possibilities to colour the first quiver. But now the second quiver gives also rise to some localization data. Dividing out the action of $S_3$ we have 6 possible colourings. Thus, in summary, we get $\chi(2+2+1,1+1+1)=24$. Note that for the last localization quiver there exists no possible colouring.\\[1ex]
Moreover, we have $\chi(3+2+0,1+1+1)=1$ when colouring the localization quiver in the middle, $\chi(3+1+1,2+1+0)=1$ when colouring the quiver on the left hand side, $\chi(2+2+1,2+1+0)=3+2=5$ when colouring the quivers on the left and right hand side. Taking into account the possibilities to obtain ordered partitions induced by the partitions, it is easy to verify $\sum_{{\bf |P_1|}=3,{\bf |P_2|}=5}\chi({\bf P_1},{\bf P_2})=204$ as predicted by Theorem \ref{kroneckerl1l2}.

\end{document}